%% file: main.tex
\documentclass[authoryear,review, usenames,dvipsnames]{elsarticle}
\usepackage[utf8]{inputenc}

\usepackage{dsfont}
\usepackage{latexsym}
\usepackage{amssymb}
\usepackage{amsmath}
\usepackage{amsthm}
\usepackage{algorithmic}
\usepackage{paralist}
\usepackage{algorithm}
\usepackage{upgreek}
\usepackage{array}
\usepackage{soul}

\usepackage{caption}
\usepackage{rotating}
\usepackage{tikz}
\usetikzlibrary{positioning}
\usepackage{blkarray}
\usepackage{mathtools}
\usepackage{pifont}

\usepackage{longtable}
\usepackage{multirow} 
\usepackage{multicol} 

\usepackage{graphicx}
\usepackage{pgfplots}
\pgfplotsset{compat=1.14}
\usepackage{multirow}

\usepackage{subcaption}

\usepackage{anysize}
\usepackage{fancyhdr}
\usepackage{color}
\usepackage{tocbibind}

\usepackage{epsfig}

\usepackage{xcolor}
\usepackage{lineno,hyperref}
\hypersetup{colorlinks = true, linkcolor = blue, anchorcolor =red, citecolor = blue,filecolor = red, urlcolor = red, pdfauthor=author}

\usepackage{rotating}
\usepackage{url}
\usepackage{lscape}
\usepackage[normalem]{ulem}

\newtheorem{assumption}{Assumption}

\newtheorem{prop}{Proposition}

\newtheorem*{lemma*}{Lemma}

\usepackage{graphicx}
\usepackage{multirow}

\usepackage{xcolor}
\usepackage{lineno,hyperref}
\hypersetup{colorlinks = true, linkcolor = blue, anchorcolor =red, citecolor = blue,filecolor = red, urlcolor = red, pdfauthor=author}
\usepackage{tikz}
\usepackage{rotating}
\usepackage{url}
\usepackage{ulem}
\usepackage{cancel}
\usepackage{soul}
\usepackage{footnote}

\numberwithin{equation}{section}
\usepackage[hang,flushmargin]{footmisc}

\author[ad000]{Natividad González-Blanco\corref{cor1}}
\ead{ngonzalez@uloyola.es}
\author[ad1]{Antonio J. Lozano}
\ead{antonio.lozano@dmat.uhu.es}
\author[ad2]{Vladimir Marianov}
\ead{marianov@ing.puc.cl}
\author[ad0,ad00]{Juan A. Mesa}
\ead{jmesa@us.es}

\address[ad000]{Department of Quantitative Methods, Universidad Loyola, Dos Hermanas, Spain}
\address[ad1]{Department of Integrated Sciences, Universidad de Huelva, Spain}
\address[ad2]{Pontificia Universidad Católica de Chile and Instituto Sistemas Complejos de Ingeniería (ISCI)}
\address[ad0]{Department of Applied Mathematics II, Universidad de Sevilla, Spain} 
\address[ad00]{IMUS, Sevilla, Spain}

\cortext[cor1]{Corresponding author}

\journal{Transportation Research Part B: Methodological}

\begin{document}

\begin{frontmatter}

\title{An Approach to the Joint Rapid and Slow Transit Network Design Problem}

\begin{abstract}

The increase in congestion in surface traffic, airborne pollution, and other environmental issues have motivated the transit authorities to promote public transit worldwide. In big cities and large metropolitan areas, adding new rapid transit lines attracts more commuters to the public system, as they frequently allow saving travel time as compared to the private mode (car) that faces high congestion. In addition, the travel time has less variability with respect to preset schedules, and rapid lines are more efficient than slow modes operated by buses. When a new rapid transit line is constructed, it partially replaces the traffic of existing slow transit lines. As a consequence, some of the slow-mode lines have to be either canceled or their routes modified to collaborate properly with the new rapid transit line. This process is usually carried out in a sequential way, thus leading to suboptimal solutions.

In this paper, we consider an integrated model for simultaneously designing rapid and redesigning slow networks. The aim of the model is community-oriented, that is, to maximize the demand covered (or captured) by both modes. We present a mathematical programming formulation that is solved by using a specially improved Benders decomposition. For this purpose, we include a partial decomposition to speed up the computation. The computational experiments are done on a case study based on real data obtained from a survey of mobility among transportation zones in the city of Seville.
\end{abstract}

\begin{keyword}
Rapid Transit\sep 
Slow Transit\sep
Network Design\sep
Partial Benders Decomposition.
\end{keyword}
\end{frontmatter}

\input{Introduction}

\input{problem_definition}

\input{problem_formulation}

\input{algorithmic_discussion}

\input{computational_experiment}

\section{Conclusions}\label{sec:conclusions}

In this paper, we introduced and studied the problem of designing a rapid transit line while, at the same time, an existing slow transit line is redesigned, considering the existence of a private mode. The aim is to maximize the attraction of traffic from the private mode to the public mode. In addition, we show that the sequential method of solving the problem leads to suboptimal solutions. We consider the network of the city of Seville as a case study. Due to the complexity of the problem, we have developed a Partial Benders decomposition approach concerning the set of the O/D pairs. The computational results show that our proposal is competitive with the existing exact methods. Further research on this problem includes the development of a metaheuristic to improve the proposed approach even more, particularly for larger instances.

\section*{Acknowledgments}

The third author is partially supported by grants ANID FONDECYT 1200706 and ANID PIA AFB220003. The fourth author is partially supported by Grant PID2020-114594GB-C21 funded by MICIU/AEI(Spain)
/10.13039/501100011033

\input{appendix.tex}

\bibliography{biblio}

\end{document}

%% file: introduction.tex
\section{Introduction}

Public Transit Network Design has been an important area of research since the 70th of last century. Increasing awareness of efficiency, mobility, sustainability, air contamination, energy consumption, and greenhouse emissions increases, has required the application of increasingly sophisticated analytical methods for planning transit systems. These sophisticated methods and models of transit network design present important challenges due to the high computational complexity and large-scale size of real problems \cite{CanEtAl2015}. 

Public transit systems can be roughly classified into slow and rapid. The main difference between both classes is the availability of a reserved right-of-way. Those systems, such as buses and trolleybuses, that share the space of the street network with other traffic as private vehicle, have a low commercial speed. However, rapid transit systems (metro, commuter trains, etc.) do not share the space with other traffic and usually do not have level crossings, the only interaction being at stations just for pedestrians. Most cities worldwide have a slow transit network, while about 250 urban areas have a metro network. Whereas slow transit networks are selected as a sub-network from the street network, a rapid transit network is often constructed from scratch, but the cost is much higher since it includes expensive infrastructure constructions. 

Rerouting bus lines is very common when a rapid transit line starts its operation. During the last five years, about 150 new metro lines have been added to metro networks around the world, and 30 new metro systems have been inaugurated. About 180 new lines have become operating. Many other existing lines have been extended or upgraded. Moreover, numerous new modern trams, train trams, light rails, and commuter lines have also recently started their operation. In almost all cases, bus lines were (partially) doing the service before, and when a rapid transit line is put into service, some bus routes could become totally or partially useless or, at least, require some redesign. One typical example is the adaptation of the Bus Rapid Transit TranSantiago when Metro Line 7 will start its operation. Another example is bus lines 1 and 3 of TUSSAM (Municipal Bus Company) with planned  Line 3 of the Seville Metro. Usually, the metro planning projects do not take into account the bus system because they often depend on different agencies. After the introduction of the rapid transit service, the bus system is reorganized. However, this procedure could lead to sub-optimal solutions. 

The agents involved in the transit network design can be classified into three groups: construction and operation companies, users, and the community in general, which is represented by the transportation agencies and authorities. Each of these groups has different objectives: cost and/or profit for the companies involved in construction and operation, travel time, price, comfort, availability and avoidance of multiple transfers, and similar other features for users, and the general interest that can be quantified by the trip coverage and/or reduction of private traffic, for the transportation agencies.
We optimize the general interest of the community,  by maximizing the joint trip coverage of both systems, which has been considered the best measure of the future expected ridership of the integrated public system.  In the network design phase of the sequential transportation planning procedure capture of passengers for this public system is often assumed to be achieved when the time to travel between an origin and a destination is less than the travel time it would take using private cars.

Transport systems can be represented as networks with nodes (stations or stops) and edges (stretches/sections). Mathematical programming programs use design variables to decide the nodes and edges to be selected/constructed, and flow variables to route the travel demand along the edges/arcs. This structure suggests the application of decomposition methods for solving problems when the size of the instances is not too large but heuristics, metaheuristics of mateheuristics for large instances.

\subsection{Literature review}
In this subsection, we will revise the research done so far on the topics involved in the problem we address in this paper: location of rapid transit alignments, rapid transit network design, transit network design, feeder-bus network design, and some related problems.

Transit Network Design has been considered as the first step in the sequential procedure for the planning of transit systems, and it is a crucial step for the whole planning. The following steps, frequency setting, timetabling, rolling-stock and crew management, etc., strongly depend on the layout of lines previously decided. Hundreds of papers are dedicated to the transit network design problem, but most deal with the bus mode. Therefore, we will refer the reader to recent literature reviews, and we only mention some papers dealing with related problems below. In \cite{GuiEtAl2008}, the research done on bus network design, including scheduling, is reviewed. The paper by \cite{FaEtAl2013} is devoted to a review of the literature on both road and bus network design problems. The topic reviewed in the paper by \cite{CanEtAl2015} is the bus line planning from the point of view of the Mathematical Programming programs. The chapter \cite{MauEtAl2021} contains an extensive discussion on the characteristics of the different public modes, several mathematical programming programs for both the physical and route network design problems, and exact and heuristics procedures for solving these programs. 

The problem of locating a rapid transit alignment consists of choosing, in a given candidate space, a set of points to construct stations and the connections (tracks) between pairs of them, so that one or several objective functions are optimized among all the feasible solutions forming a path. The first paper dealing with this problem was \cite{GenEtAl1995}, where a revision of the main criteria used both in research and technical documents is provided, and a tabu search procedure is proposed for the location of an alignment maximizing the population covered. In \cite{DuEtAl1996}, the population covered by each station is computed by assigning non-increasing weights to the population living in Manhattan-distance catchment areas around it. A complete computational experience for the tabu search procedure is described. The paper of \cite{BruEtAl1998} considers two possible networks for commuters: the bimodal pedestrian-public and the private. Each demand, given by an origin and destination, is assigned to one of these networks based on the best user´s travel cost. Two criteria are considered: total travel cost and construction cost. Furthermore, non-inferior solutions based on the k-shortest path algorithm are provided. An algorithm that consists of two phases, construction and improvement, is designed in the paper by \cite{BruEtAl2002}. This algorithm improves the running computational time of that given in \cite{DuEtAl1996}. A community-oriented objective, the total trip coverage, is used in the paper by \cite{LaEtAl2005}. The trip coverage provided by a pair of stations is the expected number of trips between both stations, and is computed by composing the attraction of each station with the share given by a mode choice model based on a binomial logit function. The proposed algorithm is a greedy approximation scheme that consists of a construction phase followed by a post-optimization procedure. A totally different criteria is applied in \cite{LaEtAl2009} in which, with the help of a modified Voronoi diagram, good alignments are generated to reduce the proximity to historical buildings. With the aim of avoiding the common zig-zag phenomenon when maximizing coverage, \cite{MarEtAl2019} introduce and minimize a measure of what they call discrete curvature.

Models and algorithms for the Rapid Transit Network Design problem have been recently reviewed in \cite{LaMe2020}. In the general setting, these problems are expected to be NP-hard. In particular, the problem of maximizing the estimated ridership subject to a construction cost, and those with the center and median objective functions, have been proved to belong to the NP-hard class (\cite{PeEtAl2020}, \cite{BuEtAl2023}). For this reason, the existing research has been oriented to the use of decomposition approaches, as Column Generation and Benders procedures, to solve problems with moderate size \cite{bucarey2022benders}, to apply metaheuristics for finding approximate solutions (\cite{CanEtAl2017}, \cite{ChenEtal2024}), or applying methodological reductions of the problems as the restriction on the lines to be in previously stated corridors (\cite{LaPas2015}, \cite{GuEtAl2018}).

The Feeder Bus Network Design (FBND) problem arises when a rapid transit system is functioning and a set of bus routes (and possibly their frequencies) must be decided to carry riders to the rapid transit stations. The feeder bus problem alone has been classified as NP-hard in \cite{MarVaz1998}. Due to its computational complexity, researchers have applied metaheuristics to solve the related problem. In \cite{AlEtAl2015}, a combination of evolutionary techniques has been applied to the problem of designing a set of feeder bus routes and determining the frequencies, aiming at minimizing the sum of operator, users and social costs. In the paper by \cite{LinWang2014}, a multiobjective approach that takes into account route lengths, travel time, and trip coverage is proposed and solved. The feeder bus planning problem has been researched to some extent \cite{deng2013optimal} even including the cost and time of riders of a fixed rapid transit system. Still, as far as the authors are aware, for the simultaneous and cooperative slow and rapid transit network design problem, no research has been done except for continuous models (\cite{FanEtAl2018}). This paper considers an idealized city with a grid street pattern with a uniform demand pattern. The objective function is the sum of both patron and agency average costs. A numerical algorithm solves the problem, and some insights are derived. 
With the purpose of filling the gap in discrete models, in this paper, an integer mathematical programming program for the integrated rapid and slow transit network design is presented, which designs the rapid transit line and redesigns the routes, if necessary, of some of the slow lines.

\subsection{Contributions of this paper}
 
The contributions of this work are the following:
 \begin{itemize}
     \item We present a mathematical formulation to address the complex task of planning the development of a rapid transit line while simultaneously relocating an existing slow transit line. This intricate process is designed to carefully consider and balance various factors, with a primary focus on optimizing the overall demand served by both transit systems. By integrating these two tasks within a unified mathematical framework, we aim to achieve a harmonious coexistence of the transit lines, ensuring that the overall transportation network is strategically enhanced, against performing sequentially both designs.

     \item Given the complex nature of the presented formulation, we have also contributed by adapting and implementing the known Benders decomposition methodology (\cite{benders1962partitioning}) to address the complexity inherent in our approach effectively. Our adaptation involves an implementation, which is currently known as the Branch and Benders cut technique.
     \item Furthermore, to refine and mitigate the time-consuming and unstable aspects associated with the aforementioned Benders decomposition procedure, we have drawn insights from \cite{belieres2020benders} and \cite{rahmaniani2017benders}. 
     
     \item Finally, the proposed procedure facilitates the solution of real-world scenarios, such as the one examined in this research, consisting in the network in the city of Seville, which is not possible to solve by the direct use of a solver (\texttt{CPLEX} in our case). As happens in real-life problems, due to the large scale of the Seville City network, we have considered different sub-instances of it to solve the corresponding problems.
 \end{itemize}

\subsection{Structure of this paper}
The structure of this paper is as follows. Section \ref{sec:problem_definition} presents the problem of designing a rapid transit line while simultaneously relocating an existing slow transit line. Furthermore, in Section \ref{sec:problem_formulation}, we propose a mathematical formulation for the previously mentioned problem. Subsequently, in Section \ref{sec:algorithmic_discussion}, given the computational complexity of the problem, we describe the development of a Benders decomposition approach, encompassing the provision of facet-defining cuts. Section \ref{sec:comp_exp} comprehensively examines a real case study of the Seville city network. Besides, this section shows an extensive computational experience related to our Benders decomposition approach proposed in the previous section. Finally, our conclusions are presented in Section \ref{sec:conclusions}.

%% file: problem_definition.tex
\section{Problem description and assumptions}\label{sec:problem_definition}

In this section, we present the problem in detail. First, we describe the problem, introduce the elements and parameters needed, and then we describe the assumptions.

\subsection{Problem description}
In order to describe the problem, we need to define the following elements and parameters.

\begin{enumerate}
\item We consider a corridor embedded in an urban area that is tessellated into transportation zones, with their corresponding centroid, being $N_c$ the whole centroid set. Given that the problem dealt with in this paper has a strategic character, we only consider one demand matrix without distinguishing the period of time, purposes, or means of transportation.
 
 \item The relevant networks to be considered are described as follows.
 
 \begin{enumerate}
 \item The street network $\mathcal{N}_{\mathcal{S}}=(N_{\mathcal{S}},E_\mathcal{S})$
 consists of those streets that can be traversed by the slow mode of transportation.
 
 \item The potential network for the rapid transit line $\mathcal{N}_{\mathcal{R}}=(N_\mathcal{R}, E_\mathcal{R})$.

 \item Since we are using centroids that concentrate the demand between pairs of transportation zones, we consider a network $\mathcal{N_W}=(N_\mathcal{W}, E_\mathcal{W}),$ which edges connect centroids with the closest candidate nodes for both of the rapid and slow networks. These edges represent the walking legs of the riders routes between the centroids and the boarding/alighting stations/stops of the rapid/slow transit network. Without loss of generality, the walking mode can be replaced by different micro-mobility modes, such as bicycle and electric scooter,  park\&ride, or aggregate them \cite{LiuOu2021}, \cite{WangEtAl2022}.
 \end{enumerate} 

Note that $\mathcal{R}$ and $\mathcal{S}$ have some nodes in common.
 
\item For modes of transportation rapid, $\mathcal{R}$, slow, $\mathcal{S},$ and walking, $\mathcal{W}$, we define $A_\mathcal{R},$ $A_\mathcal{S}$, and $A_\mathcal{W}$ as their respective sets of arcs.

\item For the rapid transit line $\mathcal{R}$, there exists a maximum number of edges $E_{\mathcal{R}}^{max}$ to be constructed. To keep within limits the disturbance to current users of the slow line, $\mathcal{S}$, upper bounds $E_{\mathcal{S}}^{max}>0$ and $E_{\mathcal{S}}^{id}\geq 0$ are given to the number of edges that will belong to the modified slow line $\mathcal{S}$, and the minimum number of edges that must remain unmodified in line $\mathcal{S}$, respectively. For that, vector $b_e^{\mathcal{S}},\, e\in E_{\mathcal{S}}$ denotes the current path of the slow line $\mathcal{S}$. 
	
\item We assume that there is a set $O_\mathcal{R}$ of possible starting points and $D_\mathcal{R}$ of possible endpoints of the rapid line. Similarly, there is a set $O_\mathcal{S}$ of possible starting points and $D_\mathcal{S}$ of possible endpoints of the slow line.

\item The set of demands $W$ is a subset of $N_c\times N_c$. The mobility pattern is given by a matrix $G = (g^w)$, where $g^w, w = (w^s,w^t)\in W$, denotes the number of expected passengers going from $w^s$ to $w^t$. 

\item For each $w\in W$, there exists a fixed cost of going from node $w^s$ to node $w^t$ using the private mode of transportation, denoted by $u^w$. 
 
 \item Let $\delta(i)$ be the set of edges of $E_{\mathcal{R}}$ incident to node $i$. The notation $\delta_+(i)$ ($\delta_-(i)$ respectively) is used to denote the set of arcs of $A_{\mathcal{R}}$ going out (in, respectively) of node $i\in N$. Similarly, we use the notation $\vartheta(k)$ and $\vartheta_+(k)$ ($\vartheta_-(k)$ respectively) to denote the set of edges of $E_{\mathcal{S}}$ incident to node $k$ and the set of arcs going out (in, respectively) of node $k$ in $A_{\mathcal{S}}$. Note that the sets of edges are related to the network design, and the sets of arcs refer to the flow paths in the designed network.
 
 \item The set of potential transfer nodes is denoted by $N_{trans}\subseteq N_{\mathcal{R}}\cap N_{\mathcal{S}}$.
	
 \item Other costs are those for walking to/from a potential location of a station/stop of one of the two modes of transportation, $t_{w^sk}$ and $t_{kw^t}$; the user cost of traversing arc $a$ in the rapid and slow mode, $t^{\mathcal{R}}_a$ and $t^{\mathcal{S}}_a$, respectively; the transfer cost at station $k$ from $\mathcal{S}$ to $\mathcal{R}$ and from $\mathcal{R}$ to $\mathcal{S}$, $t^{\mathcal{SR}}_k$ and $t^{\mathcal{RS}}_k$, respectively; the dwell (stop time) costs $t^{\mathcal{R}}_{stop}$ and $t^{\mathcal{S}}_{stop}$, which will be assumed independent from nodes since we are in the strategic phase; and the waiting time at stations/stops, $t_{wait}$, which is usually set as half of the headway. 
\end{enumerate}

Considering all of these elements and the parameters, the goal of the problem is to maximize the demand covered by the rapid and slow transit lines cooperating with each other in competition with the private mode.

\subsection{Assumptions}
\begin{assumption}
The demand is assigned in an all-or-nothing way to the private mode or to the public mode composed of the rapid transit line $\mathcal{R}$ and the slow transit line $\mathcal{S}$.
\end{assumption}

\begin{assumption}
The interspace-station distance of $\mathcal{R}$ must be at least $C_1$. However, we do not consider an inter-spacing distance for the slow line because the location of stops for bus lines depends on the available space in the street network.
\end{assumption}

\begin{assumption}
The maximum walking distances from an origin centroid to a station/stop or from a station/stop to a destination centroid are $C_2$ or $C_3$, depending on whether the station belongs to the rapid or slow mode, respectively.
\end{assumption}

\begin{assumption}
As happens in practice, usually no centroid coincides with a potential station.
\end{assumption}

\begin{assumption}
For any trip, only one transfer from slow to rapid mode and from rapid to slow is allowed. Hence, two transfers are possible for each demand.
\end{assumption}

\begin{assumption}

We only consider routes from/to centroids to stations/stops for the pedestrian mode.
\end{assumption}

\begin{assumption}
The utility of a route between each pair of centroids by the public mode is composed of the walking time from the origin centroid to a station/stop of the rapid or slow mode, the waiting time for the vehicle, the travel time using the rapid or slow mode or both, the transfer times, and the walking time from a station/stop to the destination centroid. Walking between stations/stops is not considered. 
\end{assumption}

%% file: problem_formulation.tex
\section{Problem formulation}\label{sec:problem_formulation}

In this section, we propose a formulation for the problem described in the previous Section. 

\begin{figure}[H]
\centering
\begin{subfigure}{0.5\textwidth}
        \centering
        \includegraphics[width=\linewidth]{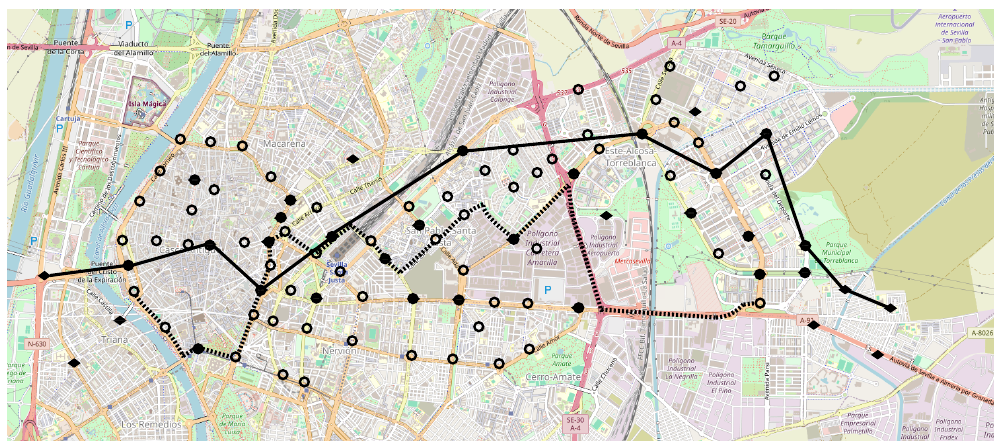}
\caption{The objective value of the sequential design is $\texttt{obj\_v}=831$.}
\end{subfigure}
 \begin{subfigure}{0.5\textwidth}
        \centering
        \includegraphics[width=\linewidth]{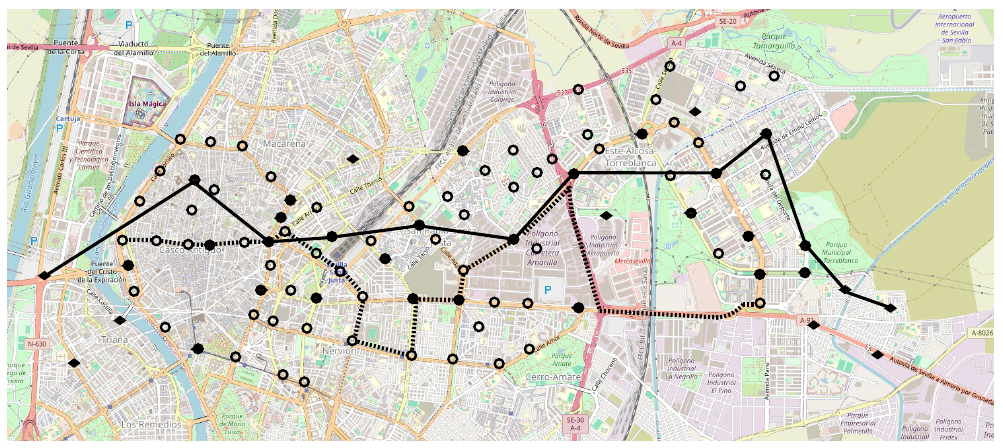}
\caption{The objective value of the integrated design is $\texttt{obj\_v}=1031$.}
\end{subfigure}
\caption{Comparison of the Sequential (above) and Integrated (below) designs for the Seville City network considering the \texttt{Instance$_{g\geq 150}$}. Continuous and dashed lines represent the designs of $\mathcal{R}$ and $\mathcal{S}$, respectively. Black nodes correspond with the set of potential stations of $\mathcal{R}$ and the white ones with the potential stations of $\mathcal{S}$.}
\label{fig:Sequential_Seville_instance}
\end{figure}

Locating each line independently without taking into account the influence that may exist between them, or even sequentially, which is the usual method in practice, can lead to suboptimal solutions. Currently, the rapid transit line $\mathcal{R}$ is located first, and then the slow line $\mathcal{S}$ is redesigned. The integrated model presented in this section results in an optimum design concerning the maximization of the coverage for the whole public transport (composed of the rapid and slow modes), as shown in Figure \ref{fig:Sequential_Seville_instance} for a sub-instance of the Seville City network (see details in Section \ref{sec:comp_exp}). The Figure shows that these different approaches result in different network designs. In what follows, we name such formulations as Sequential Network Design and Integrated Network Design (IND).

\subsection{Variables}

For the formulation proposed, we used the following set of variables:

\begin{enumerate}
	\item $x^\mathcal{R}_e = 1 \text{ if edge } e = \{k, l\} \in E_\mathcal{R}$ is included in the rapid transit line $\mathcal{R}$; $0$ otherwise. Analogously, $x^\mathcal{S}_e = 1 \text{ if edge } e = \{k, l\} \in E_\mathcal{S}$ is included in the slow transit line $\mathcal{S}$; $0$ otherwise.
	\item $y^\mathcal{R}_i = 1 \text{ if node } i\in N_{\mathcal{R}}$ is included in the alignment of the rapid system $\mathcal{R}$, but it does not stop on it; $0$ otherwise. These variables allow the inclusion of non-stop nodes in the rapid line, providing more flexibility to the design and resulting in a better relationship between speed and trip coverage.
	\item $z^\mathcal{R}_i = 1$ if $\mathcal{R}$ stops at $i$; $0$ otherwise. Analogously, $z^\mathcal{S}_k = 1$ if $k$ is a stop of mode $\mathcal{S}$; $0$ otherwise. 
 \item  $f^w = 1$ if demand $w$ uses $\mathcal{S}$, $\mathcal{R}$, or the combined modes $\mathcal{RS}$ and $\mathcal{SR}$.
	\item $f^{w\mathcal{R}}_a = 1$ if demand $w$ traverses arc $a\in A_\mathcal{R}$; $0$ otherwise. Analogously, $f^{w\mathcal{S}}_a = 1$ if demand $w$ traverses arc $a\in A_\mathcal{S}$; $0$ otherwise.
	\item $f^{w\mathcal{SR}}_k = 1$ if demand $w$ transfers from $\mathcal{S}$ to $\mathcal{R}$ at node $k\in N_{trans}$; $0$ otherwise. Analogously, $f^{w\mathcal{RS}}_k = 1$ if demand $w$ transfers from $\mathcal{R}$ to $\mathcal{S}$ at node $k\in N_{trans}$; $0$ otherwise.
        \item $v_{{w^s}k}^{w\mathcal{R}} = 1$ if demand originated at $w^s$ walks to a rapid line station located at $k$; $0$ otherwise. Analogously, $v_{{w^s}k}^{w\mathcal{S}} = 1$ if demand originated at $w^s$ walks to a slow line station located at $k$; $0$ otherwise. In the same way, we have constructed variables $v_{{kw^t}}^{w\mathcal{R}} = 1$ so state that if demand leaving the rapid line at a station located at $k$, walks to the destination $w^t$; $0$ otherwise. Analogously, $v_{{kw^t}}^{w\mathcal{S}} = 1$ if demand leaving the slow line at a stop located at $k$, walks to the destination $w^t$; $0$ otherwise. 
\end{enumerate}

\subsection{Objective function and constraints}

The problem aims to design line $\mathcal{R}$ and to re-design line $\mathcal{S}$ to maximize the joint traffic coverage or capture of public modes, thus minimizing the private traffic:
\begin{flalign}
&\max_{\boldsymbol{x},\boldsymbol{y},\boldsymbol{z},\boldsymbol{f},\boldsymbol{v}^{\mathcal{R}},\boldsymbol{v}^{\mathcal{S}}} \sum_{w \in W} g^w f^w.
\end{flalign}

\begin{itemize}
	\item Budget constraints. We approximate the constraints on budget by imposing upper bounds on the number of edges in each line. This makes the constraints tighter and increases the tractability of the problem. 
	\begin{flalign}
		& \sum_{e \in E_{\mathcal{R}}} x^{\mathcal{R}}_e \le E_{\mathcal{R}}^{max},\label{eq:budget} \\
		& \sum_{e \in E_{\mathcal{S}}} x^{\mathcal{S}}_e \le E_{\mathcal{S}}^{max}.
	\end{flalign}
	\item Design constraints. 
	\begin{flalign}
		& x^{\mathcal{R}}_e \leq z^{\mathcal{R}}_i + y^{\mathcal{R}}_i, \quad e\in E_{\mathcal{R}},\, i\in e,& \label{eq:design_1}\\
           &\sum_{o\in O_{\mathcal{R}}} z^{\mathcal{R}}_o = 1, \label{eq:design_2}\\
		&\sum_{d\in D_{\mathcal{R}}} z^{\mathcal{R}}_d = 1, \label{eq:design_3}\\           
		 &\sum_{o\in O_{\mathcal{R}}}\sum_{e\in \delta(o)} x^{\mathcal{R}}_e = 1, \label{eq:design_4}\\
		 &\sum_{d\in D_{\mathcal{R}}}\sum_{e\in \delta(d)} x^{\mathcal{R}}_e = 1, \label{eq:design_5}\\	
		& z^{\mathcal{R}}_i+y^{\mathcal{R}}_i \leq 1, \quad i\in N_{\mathcal{R}}, & \label{eq:design_6}\\
		&\sum_{e\in E_{\mathcal{R}}}x^{\mathcal{R}}_e + 1 = \sum_{i\in N_{\mathcal{R}}}(y^{\mathcal{R}}_i + z^{\mathcal{R}}_i),  & \label{eq:design_7}\\
		& \sum_{e\in\delta(k)} x^{\mathcal{R}}_e = 2(z^{\mathcal{R}}_k + y^{\mathcal{R}}_k), \quad k\in N_{\mathcal{R}}\setminus(O_{\mathcal{R}}\cup D_{\mathcal{R}}),& \label{eq:design_8} \\
		& x^{\mathcal{S}}_e \leq z^{\mathcal{S}}_i, \quad e\in E_{\mathcal{S}},\, i\in e,& \label{eq:design_9}\\
   &\sum_{o\in O_{\mathcal{S}}} z^{\mathcal{S}}_o = 1, & \label{eq:design_10}\\
		&\sum_{d\in D_{\mathcal{S}}} z^{\mathcal{S}}_d = 1, & \label{eq:design_11}\\
   &\sum_{o\in O_{\mathcal{S}}}\sum_{e\in \vartheta(o)} x^{\mathcal{S}}_e = 1, \label{eq:design_12}\\
		 &\sum_{d\in D_{\mathcal{S}}}\sum_{e\in \vartheta(d)} x^{\mathcal{S}}_e = 1, \label{eq:design_13}
   \end{flalign}
    \begin{flalign}
		&\sum_{e\in E_{\mathcal{S}}}x^{\mathcal{S}}_e + 1 = \sum_{i\in N_{\mathcal{S}}}z^{\mathcal{S}}_i,& \label{eq:design_14}\\
		& \sum_{e\in\vartheta(k)} x^{\mathcal{S}}_e = 2z^{\mathcal{S}}_k, \quad k\in N_{\mathcal{S}}\setminus(O_{\mathcal{S}}\cup D_{\mathcal{S}}), & \label{eq:design_15} \\	
  &  \sum_{e\in E_s} b_e^{\mathcal{S}} x_e^{\mathcal{S}} \ge E_{\mathcal{S}}^{id}.\label{eq:design_16}
	\end{flalign}

 Constraints \eqref{eq:design_1} force the extremes of a constructed edge of the line ${\mathcal{R}}$ to be stations (constructed nodes) or non-stop nodes. Constraints \eqref{eq:design_2} and \eqref{eq:design_3} impose that exactly one node has to be selected from the sets of possible origins and destinations of the rapid transit line. Constraints \eqref{eq:design_4} and \eqref{eq:design_5} enforce just one incident edge to $O_R$ and $D_R$ is constructed. 
  Constraints \eqref{eq:design_6} do not allow a node $i\in N_{\mathcal{R}}$ to be simultaneously a stop and non-stop node. Constraints \eqref{eq:design_7} impose that the rapid line must be a forest graph, and \eqref{eq:design_8} that the degree of non-terminal nodes is 2. These constraints along with \eqref{eq:design_2}, \eqref{eq:design_3}, \eqref{eq:design_4} and \eqref{eq:design_5} preclude the existence of cycles. Thus, the forest tree has just one component, which is a chain graph.  We consider similar constraints for the re-location of the slow line $\mathcal{S}$. Constraints \eqref{eq:design_9} enforce the extremes of a constructed edge of line $\mathcal{S}$ to be stop nodes. Constraints \eqref{eq:design_10} to \eqref{eq:design_15} impose the slow line to be represented as a chain-graph.
  Constraint \eqref{eq:design_16} states that the old and new line $\mathcal{S}$ must coincide in a given number of edges. That is, changes in the bus routes can be made, but the new path can coincide partially with the old one.
	
	\item Relation between mode choice and pedestrian choice variables.
	\begin{flalign}
            & \sum\limits_{k\in N_{\mathcal{R}}} v_{{w^s}k}^{w\mathcal{R}} +\sum\limits_{k\in N_{\mathcal{S}}} v_{{w^s}k}^{w\mathcal{S}}= f^w, \quad w\in W, \label{eq:relation_3}\\
            & \sum\limits_{k\in N_{\mathcal{R}}} v_{k{w^t}}^{w\mathcal{R}} +\sum\limits_{k\in N_{\mathcal{S}}} v_{k{w^t}}^{w\mathcal{S}}= f^w, \quad w\in W,\label{eq:relation_4}\\
 		& v_{{w^s}k}^{w\mathcal{R}}\leq z^{\mathcal{R}}_{k},\quad w\in W, \quad  k\in N_{\mathcal{R}} ,& \label{eq:relation_5}\\
		& v_{{w^s}k}^{w\mathcal{S}}\leq z^{\mathcal{S}}_{k},\quad w\in W, \quad  k\in N_{\mathcal{S}} ,& \label{eq:relation_6}\\
		& v_{k{w^t}}^{w\mathcal{R}}\leq z^{\mathcal{R}}_{k},\quad w\in W, \quad  k\in N_{\mathcal{R}} ,& \label{eq:relation_7}\\
		& v_{k{w^t}}^{w\mathcal{S}}\leq z^{\mathcal{S}}_{k},\quad w\in W, \quad  k\in N_{\mathcal{S}}. &\label{eq:relation_8}
  \end{flalign}
 
  Constraints \eqref{eq:relation_3} state that each outflow from $w^t$  must use the walking mode to a location $k$ at most belonging to one of the two public modes of transport. In the same way, constraints \eqref{eq:relation_4} state that each inflow of $w^s$ must use the walking mode from a location $k$ at most belonging to one of the two public modes of transport. Constraints \eqref{eq:relation_5} and \eqref{eq:relation_6} impose that if a station $k$ is not a stop, then the outflow of node $w^s$ walking to it cannot be satisfied. In the same way, constraints \eqref{eq:relation_7} and \eqref{eq:relation_8} enforce that if at location $k$ there is not a stop node, then the inflow from it to $w^t$ cannot be satisfied. 
  
\item Flow conservation constraints in non-transfer stations. Flows have to be maintained either by slow or rapid modes.
 \begin{flalign}
  & \sum_{a \in \delta_w^+(k)} f^{w\mathcal{R}}_a - \sum_{a \in \delta_w^-(k)} f^{w\mathcal{R}}_a = v_{{w^s}k}^{w\mathcal{R}}-v_{k{w^t}}^{w\mathcal{R}}, \quad w\in W, \quad k\in N_{\mathcal{R}}\setminus N_{trans},\label{eq:flow_1}\\
  & \sum_{a \in \vartheta_w^+(k)} f^{w\mathcal{S}}_a - \sum_{a \in \vartheta_w^-(k)} f^{w\mathcal{S}}_a = v_{{w^s}k}^{w\mathcal{S}}-v_{k{w^t}}^{w\mathcal{S}}, \quad w\in W, \quad k\in N_{\mathcal{S}}\setminus N_{trans}.\label{eq:flow_11}\\
					& v_{{w^s}k}^{w\mathcal{R}} \leq \sum_{a \in \delta_w^+(k)} f^{w\mathcal{R}}_a,\quad w\in W, \, k\in N_{\mathcal{R}},& \label{eq:flow_2}\\
					& v_{{w^s}k}^{w\mathcal{S}} \leq \sum_{a \in \vartheta_w^+(k)} f^{w\mathcal{S}}_a,\quad w\in W, \, k\in N_{\mathcal{S}},& \label{eq:flow_3}\\
					& v_{k{w^t}}^{w\mathcal{R}} \leq \sum_{a \in \delta_w^-(k)} f^{w\mathcal{R}}_a,\quad w\in W, \, k\in N_{\mathcal{R}}, \label{eq:flow_4}&\\
					& v_{k{w^t}}^{w\mathcal{S}} \leq \sum_{a \in \vartheta_w^-(k)} f^{w\mathcal{S}}_a,\quad w\in W, \, k\in N_{\mathcal{S}}, &\label{eq:flow_5}\\
      &f^{w\mathcal{R}}_a\leq f^w, \quad w\in W, a\in A_{\mathcal{R}},& \label{eq:flow_6}\\
		&f^{w\mathcal{S}}_a\leq f^w, \quad w\in W, a\in A_{\mathcal{S}}.&\label{eq:flow_7}
  \end{flalign}

  Constraints \eqref{eq:flow_1} and \eqref{eq:flow_11} impose that the difference between the outflow and inflow of $w$ at node $k,$ either by the rapid or slow mode must be equal to the corresponding difference of the pedestrian mode at station/stop $k.$
  Constraints \eqref{eq:flow_2} and \eqref{eq:flow_3} guarantee that the flow outcoming from location $w^s$ to station $k$ only can use one arc going out from $k$. In the same way, constraints \eqref{eq:flow_4} and \eqref{eq:flow_5} guarantee that the flow incoming from station $k$ to location $w^t$ only can use one arc going into $k$. Constraints \eqref{eq:flow_6} and \eqref{eq:flow_7} impose that if the demand $w$ is not captured, then no flow traverses any edge of any of the public modes. 

	\item Transfer constraints. Only one transfer from slow to rapid mode and from rapid to slow is allowed, by using constraints \eqref{eq:transfer_1} and \eqref{eq:transfer_2}. Besides, flow conservation constraints for transfer points \eqref{eq:transfer_3}, \eqref{eq:transfer_4}, \eqref{eq:transfer_5} and \eqref{eq:transfer_6} are
needed.
	\begin{flalign}
		&\sum_{k \in N_{trans}} f^{w\mathcal{SR}}_k \leq 1, \quad w\in W,  &\label{eq:transfer_1}\\
		&\sum_{k \in N_{trans}} f^{w\mathcal{RS}}_k \leq 1, \quad w\in W,  &\label{eq:transfer_2}
  \end{flalign}
  \begin{flalign}
		& \sum_{a \in \delta_w^-(k)} f^{w\mathcal{R}}_a + f^{w\mathcal{SR}}_k - \left( \sum_{a \in \delta_w^+(k)} f^{w\mathcal{R}}_a + f^{w\mathcal{RS}}_k\right) = - v_{{w^s}k}^{w\mathcal{R}}+v_{k{w^t}}^{w\mathcal{R}}, \,\,\, w\in W, \, k\in N_{trans},& \label{eq:transfer_3}\\
		&\sum_{a \in \vartheta_w^-(k)} f^{w\mathcal{S}}_a + f^{w\mathcal{\mathcal{RS}}}_k - \left(\sum_{a \in \vartheta_w^+(k)} f^{w\mathcal{S}}_a + f^{w\mathcal{SR}}_k\right) = -v_{{w^s}k}^{w\mathcal{S}} + v_{k{w^t}}^{w\mathcal{S}}, \,\,\, w\in W, \, k\in N_{trans},& \label{eq:transfer_4}
  \end{flalign}
 \begin{flalign}
& f^{w\mathcal{\mathcal{RS}}}_k \leq \sum_{a \in \vartheta_w^+(k)} f^{w\mathcal{S}}_a, \qquad f^{w\mathcal{\mathcal{RS}}}_k \leq \sum_{a \in \delta_w^-(k)} f^{w\mathcal{R}}_a, \quad w\in W, \, k\in N_{trans} \label{eq:transfer_5} &\\
& f^{w\mathcal{\mathcal{SR}}}_k \leq \sum_{a \in \vartheta_w^-(k)} f^{w\mathcal{S}}_a, \qquad f^{w\mathcal{\mathcal{SR}}}_k \leq \sum_{a \in \delta_w^+(k)} f^{w\mathcal{R}}_a, \quad w\in W, \, k\in N_{trans}. \label{eq:transfer_6} &
 \end{flalign}

	\item Location-allocation constraints. Constraints \eqref{eq:loc_allo_1} and \eqref{eq:loc_allo_2} do not allow flow on an edge of the lines
$\mathcal{R}$ and $\mathcal{S}$ unless the edge has been constructed/chosen. In addition, each edge can be used at most in one direction by each O/D pair. Furthermore, edges belonging to both modes of transport can be used only for one of the modes by each O/D pair. With these constraints, the flow is prevented from making meaningless movements.
\begin{flalign}
		&f^{w\mathcal{R}}_a + f^{w\mathcal{R}}_{\hat{a}} \leq x^{\mathcal{R}}_e, \quad w\in W, e= \{i,j\}\in E_{\mathcal{R}}: a=(i,j), \hat{a}=(j,i),& \label{eq:loc_allo_1} \\
		&f^{w\mathcal{S}}_a + f^{w\mathcal{S}}_{\hat{a}} \leq x^{\mathcal{S}}_e, \quad w\in W, e= \{i,j\}\in E_{\mathcal{S}}: a=(i,j), \hat{a}=(j,i),& \label{eq:loc_allo_2} \\
  &f^{w\mathcal{R}}_a + f^{w\mathcal{R}}_{\hat{a}}+f^{w\mathcal{S}}_a + f^{w\mathcal{S}}_{\hat{a}} \leq 1, \quad w\in W, e= \{i,j\}\in E_{\mathcal{R}}\cap E_{\mathcal{S}}: a=(i,j), \hat{a}=(j,i).& \label{eq:loc_allo_3} 
	\end{flalign}
        \item Transfer at stations constraints. Constraints \eqref{eq:alignment_1} and \eqref{eq:alignment_2} preclude transfers from the rapid to the slow mode and vice-versa at a node unless there are stations of both modes at that node.
	\begin{flalign}
		&f^{w\mathcal{SR}}_k + f^{w\mathcal{RS}}_k\leq z^{\mathcal{R}}_k, \quad w\in W, k\in N_{trans},& \label{eq:alignment_1}\\
		&f^{w\mathcal{SR}}_k + f^{w\mathcal{\mathcal{RS}}}_k\leq z^{\mathcal{S}}_k, \quad w\in W, k\in N_{trans}.& \label{eq:alignment_2}
	\end{flalign}
	\item Mode choice constraints. They assign the demand either to the public mode or to the private one depending on the total time of the trip. Concerning the public mode, arrival times walking from the origin to a station, departure times walking from the arrival station to the destination, times between stops, transfer times, and stop times have been taken into account.
	\begin{flalign}
		\begin{split}
			&\hspace{-0.6cm}\sum_{k \in N_{\mathcal{R}}} \left( t_{{w^s}k}v_{{w^s}k}^{w\mathcal{R}} + t_{k{w^t}}v_{k{w^t}}^{w\mathcal{R}} \right) + \sum_{k \in N_{\mathcal{S}}} \left( t_{{w^s}k}v_{{w^s}k}^{w\mathcal{S}} + t_{k{w^t}}v_{k{w^t}}^{w\mathcal{S}} \right) + \sum_{a\in A_{\mathcal{R}}}t^{\mathcal{R}}_af^{w\mathcal{R}}_a + \\
			&\hspace{-0.6cm}+ \sum_{a\in A_{\mathcal{S}}}t^{\mathcal{S}}_af^{w\mathcal{S}}_a + \sum_{k\in N_{trans}}t^{\mathcal{RS}}_kf^{w\mathcal{\mathcal{RS}}}_k + \sum_{k\in N_{trans}}t^{\mathcal{SR}}_kf^{w\mathcal{\mathcal{SR}}}_k + t^{\mathcal{R}}_{stop}\sum_{k\in N_{\mathcal{R}}}z^{\mathcal{R}}_k\sum_{a\in\delta^+(k)}f^{w\mathcal{R}}_a +\\
			&\hspace{-0.6cm} + t^{\mathcal{S}}_{stop}\sum_{k\in N_{\mathcal{R}}}z^{\mathcal{S}}_k\sum_{a\in\vartheta^+(k)}f^{w\mathcal{S}}_a + f^w\left(t^{\mathcal{R}}_{wait}-\frac{1}{2}t^{\mathcal{R}}_{stop}\right)\leq u^w\,f^w, \quad w\in W.\label{eq:utility}
	\end{split}\end{flalign}

 Note that this constraint is quadratic. The linearization of these products of variables is detailed in \ref{append:linearization}.
 
 \item Maximum walking distance from the origins and to the destinations
 \begin{flalign}
     & t^{\mathcal{R}}_{w^sk}v^{w\mathcal{R}}_{w^sk}\leq C_2,\, w\in W,\, k\in N_{\mathcal{R}}, \qquad t^{\mathcal{R}}_{kw^t}v^{w\mathcal{R}}_{kw^t}\leq C_2,\, w\in W,\, k\in N_{\mathcal{R}},&\\
     & t^{\mathcal{S}}_{w^sk}v^{w\mathcal{S}}_{w^sk}\leq C_3, \, w\in W,\, k\in N_{\mathcal{S}}, \qquad t^{\mathcal{S}}_{kw^t}v^{w\mathcal{S}}_{kw^t}\leq C_3, \, w\in W,\, k\in N_{\mathcal{S}}.&
 \end{flalign}
 
 \item Constraint enforcing a minimum distance between stations. This constraint indicates that no pair of rapid transit stations can be closer than a preset distance from each other, as it would be inefficient due to both the longer total stop time, and the cost of mutual cannibalization of the demand of closeby stations.
	\begin{flalign}	
	&z^{\mathcal{R}}_i + \sum_{\substack{j \in N_{\mathcal{R}}: \\ d_{ij} \leq C_1}} z^{\mathcal{R}}_j \leq 1, \quad  \forall\, i \in N_{\mathcal{R}}. &
          \end{flalign}

\item Binary variables. All the variables are assumed to be in $\{0, 1\}$.
	\begin{flalign}
		& x^{\mathcal{R}}_e,\, x_e^{\mathcal{S}},\, y^{\mathcal{R}}_k,\, z^{\mathcal{R}}_k,\, z_k^{\mathcal{S}},\, f_a^{w\mathcal{R}},\, f_a^{w\mathcal{S}},\, f_k^{w\mathcal{RS}},\, f_k^{w\mathcal{SR}},\, f^w, v_{{w^s}k}^{w\mathcal{R}}, v_{{w^s}k}^{w\mathcal{S}}, v_{k{w^t}}^{w\mathcal{R}}, v_{k{w^t}}^{w\mathcal{S}} \in \{0,1\}.&
	\end{flalign}
\end{itemize}

Finally, the following \textit{rapid transit line shape} constraints can be added to avoid the zigzag line shape. They avoid routes with acute angles.
	\begin{flalign}	
	& f^{w\mathcal{R}}_a + \sum_{\substack{b \in \delta^w_+(k):\\\text{ angle a-b } \leq 90 }} f^{w\mathcal{R}}_b \leq 1, \quad \forall\,k\in N_{\mathcal{R}}, \, a \in \delta^w_-(k) &
    \end{flalign}
    In our case, since there is a budget constraint formulated in terms of the number of stations and links and the competition between the lines to be designed and an existing private mode of transport, zigzagging is very unlikely.

    The mathematical formulation of the problem does not aim to offer each demand pair its shortest path from origin to destination. Rather, it is designed to maximize the demand for which the total travel time through the public network is shorter than the travel time using the private mode.

%% file: algorithmic_discussion.tex
\section{Algorithmic discussion}\label{sec:algorithmic_discussion}

Formulation (IND) involves a large number of flow variables when the set of O/D pairs is extensive. To address this issue, we explore a stabilized Benders decomposition based on a Partial Benders decomposition. In this section, we justify that the formulation (IND) possesses the suitable separable structure needed for a Benders decomposition approach. Indeed, when the design variables, the mode choice variables and the walking centroids variables, namely $x^{\mathcal{R}}$, $x^{\mathcal{S}}$, $y^{\mathcal{R}}$, $z^{\mathcal{R}}$, $z^{\mathcal{S}}$, $f^w$, $v^{w\mathcal{R}}$, $v^{w\mathcal{S}}$ are fixed, the problem can be divided into $|W|$ sub-problems. Each of these sub-problems establishes the flow variables $f^{w\mathcal{R}}$, $f^{w\mathcal{S}}$, $f^{w\mathcal{RS}}$ and $f^{w\mathcal{SR}}$ as a linear problem by relaxing their integrality condition.

\subsection{Benders decomposition}\label{subsec:benders}

Benders decomposition has been previously analyzed and developed for covering problems in \cite{cordeau2019benders} and \cite{bucarey2022benders} for different purposes. The first of them is a location problem and the second one is in the network design field, particularly infrastructure design. We investigate an implementation of the branch-and-Benders-cut algorithm (\texttt{B\&BC}).

To generate cuts, we first need to relax the integrality condition on the flow variables. Proposition \ref{prop:projection} shows that this can be done without loss of generality. Let (IND\_R) denote the relaxed formulation of (IND) in which constraints $f^{w\mathcal{R}}$, $f^{w\mathcal{S}}$, $f^{w\mathcal{RS}}$, $f^{w\mathcal{SR}}\in\{0,1\}, w\in W$ are replaced by non-negativity constraints, i.e.
\begin{flalign}
 & f^{w\mathcal{R}}_a\geq 0, \quad w \in W, a \in A_{\mathcal{R}}, \\
 & f^{w\mathcal{S}}_a\geq 0, \quad w \in W, a \in A_{\mathcal{S}},\\
 & f^{w\mathcal{RS}}_k, f^{w\mathcal{SR}}_k \geq 0, \quad w \in W, k \in N_{trans}.
\end{flalign}

Let $Q$ be a set of points $(\boldsymbol{x},\boldsymbol{z}) \in \mathbb{R}^n \times \mathbb{R}^m$. Then, the projection of $Q$ onto the $x$-space, denoted by $Proj_x(Q)$, is the set of points given by: $Proj_{\boldsymbol{x}} (Q) = \{\boldsymbol{x} \in \mathbb{R}^m: (\boldsymbol{x},\boldsymbol{z}) \in Q \text{ for some } \boldsymbol{z} \in \mathbb{R}^n\}.$ Let us denote by $\mathcal{F}(IND)$ the set of feasible points of formulation $(IND)$.
 
\begin{prop}\label{prop:projection}
The projections of (IND) and (IND\_R) onto the $\boldsymbol{f}$-space coincide.
$$Proj_{\boldsymbol{x},\boldsymbol{y},\boldsymbol{f},\boldsymbol{v}}(\mathcal{F}(IND)) = Proj_{\boldsymbol{x},\boldsymbol{y},\boldsymbol{f},\boldsymbol{v}}(\mathcal{F}(IND\_R)).$$
\end{prop}

\begin{proof}
First, $\mathcal{F}(IND) \subseteq \mathcal{F}(IND\_R)$ implies $Proj_{\boldsymbol{x}^{\mathcal{R}},\boldsymbol{x}^{\mathcal{S}},\boldsymbol{y}^{\mathcal{R}},\boldsymbol{z}^{\mathcal{R}},\boldsymbol{z}^{\mathcal{S}},\boldsymbol{f},\boldsymbol{v}^{\mathcal{R}},\boldsymbol{v}^{\mathcal{S}}}(\mathcal{F}(IND)) \subseteq$ \newline $Proj_{\boldsymbol{x}^{\mathcal{R}},\boldsymbol{x}^{\mathcal{S}},\boldsymbol{y}^{\mathcal{R}},\boldsymbol{z}^{\mathcal{R}},\boldsymbol{z}^{\mathcal{S}},\boldsymbol{f},\boldsymbol{v}^{\mathcal{R}},\boldsymbol{v}^{\mathcal{S}}}(\mathcal{F}(IND\_R))$. Second, let $(\boldsymbol{x}^{\mathcal{R}},\boldsymbol{x}^{\mathcal{S}},\boldsymbol{y}^{\mathcal{R}},\boldsymbol{z}^{\mathcal{R}},\boldsymbol{z}^{\mathcal{S}},\boldsymbol{f},\boldsymbol{v}^{\mathcal{R}},\boldsymbol{v}^{\mathcal{S}})$ be a point belonging to $Proj_{\boldsymbol{x}^{\mathcal{R}},\boldsymbol{x}^{\mathcal{S}},\boldsymbol{y}^{\mathcal{R}},\boldsymbol{z}^{\mathcal{R}},\boldsymbol{z}^{\mathcal{S}},\boldsymbol{f},\boldsymbol{v}^{\mathcal{R}},\boldsymbol{v}^{\mathcal{S}}}(\mathcal{F}(IND\_R))$. For every O/D pair $w \in W$ such that $f^w = 0$ then $\boldsymbol{f}^{w\mathcal{R}}=\boldsymbol{f}^{w\mathcal{S}}=\boldsymbol{f}^{w\mathcal{RS}}=\boldsymbol{f}^{w\mathcal{SR}}=0$. In the case in which $f^w=1$, due to constraints \eqref{eq:flow_1}-\eqref{eq:utility}, flows $f^{w\mathcal{R}}_a$ and $f^{w\mathcal{S}}_a$ must be all equal to 1 on one path and 0 on the remaining paths. Hence, we show that $(\boldsymbol{x}^{\mathcal{R}},\boldsymbol{x}^{\mathcal{S}},\boldsymbol{y}^{\mathcal{R}},\boldsymbol{z}^{\mathcal{R}},\boldsymbol{z}^{\mathcal{S}},\boldsymbol{f},\boldsymbol{v}^{\mathcal{R}},\boldsymbol{v}^{\mathcal{S}})$ also belongs to $Proj_{\boldsymbol{x}^{\mathcal{R}},\boldsymbol{x}^{\mathcal{S}},\boldsymbol{y}^{\mathcal{R}},\boldsymbol{z}^{\mathcal{R}},\boldsymbol{z}^{\mathcal{S}},\boldsymbol{f},\boldsymbol{v}^{\mathcal{R}},\boldsymbol{v}^{\mathcal{S}}}(\mathcal{F}(IND))$. 
\end{proof}

Based on Proposition \ref{prop:projection}, we propose a Benders decomposition where flow variables are projected out from the model and replaced by dynamically generated Benders feasibility cuts.

\subsection{Improved master problem}\label{subsec:improved}

As it is explained in \cite{rahmaniani2017benders}, a straightforward application of the classical Benders decomposition may require excessive computing time and memory. The authors identify, among its main drawbacks, that it is time-consuming, generates poor cuts, ineffective initial iterations, and slow convergence. They define a four-dimension taxonomy of algorithmic improvements, that includes the following factors: decomposition strategy, solution procedure, solution generation and cut generation. We will focus on the first of them. It is known in the literature that the Benders decomposition method causes the master problem to lose all the information associated with the non-complicating variables. This results in instability and in a large number of iterations. 

With the purpose of improving this issue we include explicit information from the subproblems in the master. That is, instead of relaxing the whole set of flow variables, we keep part of the set into the master problem. In particular, a certain subset of the O/D pairs is selected, and all the flow variables associated with this subset are retained in the master. In Section \ref{subsec:comp_exp_benders}, it is shown that it is computationally beneficial. In addition,  we have tested whether it is more advantageous for this set to include O/D pairs with high demand, those with lower demand, or a random selection.

This idea was inspired by \cite{belieres2020benders} and \cite{crainic2021partial}, whose authors called this procedure Partial Benders decomposition, although they address problems that are different from our problem at hand. The first one of these papers presents and studies a Logistics Service Network Design Problem inspired by the management of restaurant supply chains. The second one applies it to a stochastic network design.

%% file: computational_experiment.tex
\section{Computacional experiments. Case study}\label{sec:comp_exp}

The purpose of this section is to assess the applicability of the improved Benders decomposition approach to solve the problem in a real context. We consider the case, based on real data, of the design of the trace of a rapid transit line in a corridor of the planned, but still in the study phase, of Line 2 of the Metro of Seville. In addition, we design the new routes for the current bus system. As a further aim, and due to the inherent complexity of the formulation and the large scale of the real problem, we evaluate the different configurations of the partial Benders decomposition proposed in Section \ref{sec:algorithmic_discussion}. This implementation is used as a sub-routine in a branch-and-Benders-cut scheme. This scheme allows cutting infeasible solutions along the branch-and-bound tree. For clarity, we will refer to this algorithmic routine as \texttt{B\&BC\_P}.

In our implementation, we used the \texttt{LazyConstraint\-Callback} function of \texttt{CPLEX} to separate integer solutions. Fractional solutions were separated using the \texttt{UserCutCallback} function.

Our comparative analysis involves assessing the performance of our \texttt{BD\_P} implementation against two benchmarks: the direct utilization of the \texttt{CPLEX} solver and the automatic Benders procedure proposed by \texttt{CPLEX}, denoted as \texttt{Auto\_BD}. 

We conduct our experiments on a computer equipped with an Intel Core i$5$-$7300$ CPU processor, with $2.50$ gigahertz, $4$-core, and $16$ gigabytes of RAM. The operating system used was the 64-bit Windows 10. The codes were implemented in Python 3.8 and executed using the \texttt{CPLEX} 12.10 solver through its Python interface. The \texttt{CPLEX} parameters were set to their default values, and the model was optimized in single-threaded mode.

For that, \texttt{t} denotes the value of the solution time in seconds, \texttt{gap} denotes the relative optimality gap in percent (the relative percent difference between the best solution and the best bound obtained within the time limit), \texttt{cuts} is the number of cuts generated by Benders decomposition and \texttt{obj} is the objective value.

To achieve the stated objectives, this section is structured as follows. In Subsection \ref{subsec:comp_exp_seville_instance}, we describe in detail the case study of the Seville city context. In Subsection \ref{subsec:com_exp_preliminary}, we show that the direct use of \texttt{CPLEX} and the \texttt{AUTO\_BD} are not able to get solutions, or those obtained are not competitive against the use of our ad-hoc implementation of a branch-and-Benders-cut decomposition, referred as \texttt{B\&BC}. Then, in Subsection \ref{subsec:comp_exp_benders}, we improve the ad-hoc implementation by using the Partial branch-and-Benders-cut decomposition approach, named \texttt{B\&BC\_P} in the following. Finally, in Subsection \ref{subsec:Comp_exp_case}, the main objective of this section is achieved by using the best configuration of our \texttt{B\&BC\_P}. We solve the problem for the Seville transportation system. Note that the purpose of Subsections \ref{subsec:com_exp_preliminary} and \ref{subsec:comp_exp_benders} having a computational character, justify the use of our Partial Benders decomposition, and get its best configuration.

\subsection{Seville instance}\label{subsec:comp_exp_seville_instance}
The planned metro network for Seville was based on a technical project commissioned by the Junta de Andalucía (see \cite{juntaandalucia2002}). The network design consists of four lines. Up to now, there is only one line operating (Line 1), and construction of Line 3 started at the beginning of 2024. For several reasons, Line 3 was prioritized over Lines 2 and 4, the latter being a circular line. Line 2 is currently under technical studies, but it is not expected to be constructed before 2028. For the technical project mentioned above, a mobility survey was carried out, where the city was divided into 164 transportation zones, and the daily trips between each pair of transportation zones were estimated. Up to now, there has been no other mobility survey useful for the purposes of designing transit lines. Though the metropolitan area has grown to 1.5 million inhabitants, the city's population has remained at around 700.000 inhabitants. For each transportation zone, we have calculated its centroid and concentrated the demand on it. We have considered corridors for metro line 2 and for bus line 27 of the municipal transit service TUSSAM, whose route must be moved once metro line 2 starts its operation. 

For the design of the underlying network, from which we have to choose the routes of both lines, we have chosen feasible points to be candidates for metro stations and stops on the bus line. Among the 97 nodes, 26 are shared between both modes of transport, leaving 36 nodes for the rapid transit mode and 87 for the slow transit mode. As for the edges, the two modes share 17 of the 247 total edges, with 92 belonging only to the rapid transit mode and 173 to the slow transit mode. Besides, we have considered 73 different centroids. The set of O/D pairs $W$ was formed by all possible ones, having in total 5256 O/D pairs.

The maximum number of edges $E_{\mathcal{R}}^{\max}$ to be constructed was fixed to 11. The maximum number of edges $E_{\mathcal{S}}^{\max}$ to compose the slow transit line $\mathcal{S}$ was set to 16, and the minimum number of edges that must remain unmodified in it was set to 2.

\begin{figure}[h!]
    \centering
    \begin{subfigure}{0.8\textwidth}
        \centering
        \includegraphics[width=\linewidth]{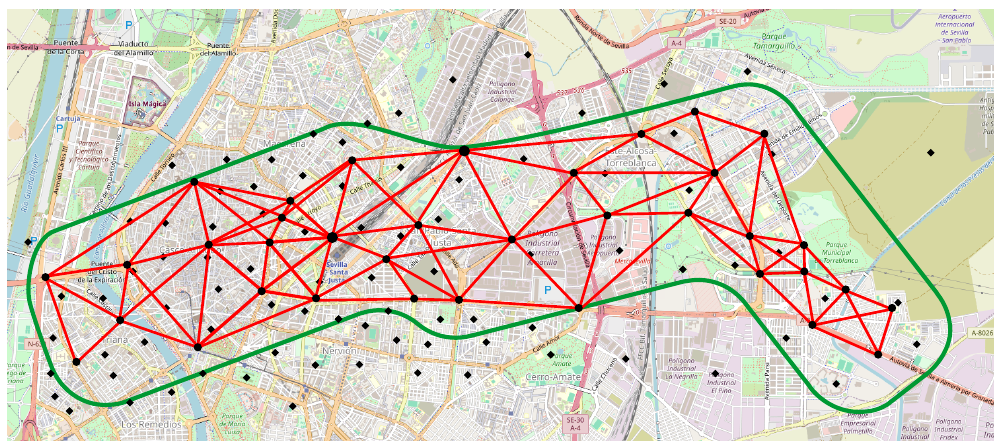}
        \caption{The potential network $\mathcal{R}$ is shown in red and the corridor in color green. The black nodes represent the centroids of transportation zones.}
        \label{fig:subfigura1}
    \end{subfigure}
    \begin{subfigure}{0.8\textwidth}
        \centering
        \includegraphics[width=\linewidth]{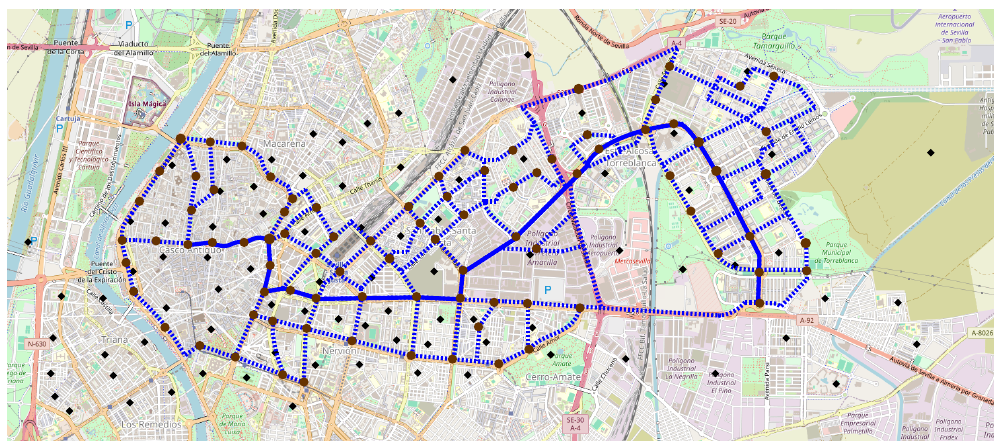}
        \caption{The potential network $\mathcal{S}$ is shown by the blue dashed line and the old path of mode $\mathcal{S}$ is depicted in blue. Besides, the black nodes represent the centroids of transportation zones.}
        \label{fig:subfigura2}
    \end{subfigure}
    \caption{Seville instance}
    \label{fig:Seville_instance}
\end{figure}

Figure \ref{fig:Seville_instance} illustrates the situation in question. On the top, the potential network $\mathcal{R}$ is shown in red and the corridor in green. The black nodes represent the sets of potential origins and destinations. That is, the centroids of transportation zones. On the bottom, the potential network $\mathcal{S}$ is represented by the blue dashed line, and the current path of mode $\mathcal{S}$ is depicted in blue. The set of centroids on the potential network is also displayed.

Each O/D pair has an associated demand with a value within the range $[0,465]$. Out of the 5256 pairs, 4114 have a demand lower than 40, with 251 of them having zero demand. These latter ones will, therefore, be excluded from the analysis.  

Parameter $C_1$ was fixed equal to 500 meters. The parameter $C_2$, which denotes the maximum distance to walk from an origin centroid to a station/stop or from a station/stop to a destination centroid in $\mathcal{R}$, was fixed to 400 meters. The parameter $C_3$, which denotes the maximum distance to walk from an origin centroid to a station/stop or from a station/stop to a destination centroid in $\mathcal{S}$, was fixed to 300 meters. Parameters $t_{w^sk}$ and $t_{kw^t}$, $w\in W, k\in N_{\mathcal{R}\cap\mathcal{S}}$, have been set considering the required time when walking at a speed of 5km/h.

Parameter $t^{\mathcal{R}}_a$, $a\in A_{\mathcal{R}}$ refers to the time to traverse arc $a$ using mode $\mathcal{R}$. It has been set to be the time to traverse arc $a$ at a speed of 70km/h. Similarly, parameter $t^{\mathcal{S}}_a$, $a\in A_{\mathcal{S}}$ has been set to be the time to traverse arc $a$ at a speed of 25km/h. 

We consider a private utility $u^w,\,w\in W$ to be equal to twice the time it takes to drive from O to D at a speed of 30 km/h.

The transfer time parameters $t^{\mathcal{RS}}_k$ and $t^{\mathcal{SR}}_k$ have been set to be equal to 9.5 and 5.5 minutes, respectively. Finally, the rest of the parameters have been set as: $t^{\mathcal{R}}_{stop}=0.5$ min, $t^{\mathcal{R}}_{stop}=1$ min and $t_{wait}=2$ min.

Due to the large scale of the instance, to test the procedures, we consider some smaller sub-instances of the full instance described previously. We build the sub-instances using sub-sets of the set of O/D pairs. That is, the set of potential nodes and edges is the same, but the O/D pairs set of each of these sub-instances corresponds to a subset of the whole potential set of O/D pairs. To construct several instances, we will first sort the O/D pairs in an ascending order given by the demand parameter $g^w$. Then, \texttt{Instance$_{g\geq150}$} is composed of the O/D pairs whose demand is equal to or bigger than 150, which is contained in the instance that considers the demand bigger than 130, and so on. To test several sub-instances, we have considered the ones in Table \ref{tab:proportion_OD_pairs}. The second row of the table refers to the number of O/D pairs that compose each sub-instance. 

\begin{table}[H]
\centering
\begin{tabular}{c|cccccc}
\cline{2-7}
& \multicolumn{6}{c}{$i$} \\
\hline
\texttt{Instance$_{g\geq i}$}   & 150 & 130 & 100 & 80 & 60 & 40 \\
number of O/D pairs &  82  &   122  & 240 & 384 & 632 & 1142 \\ 
\hline
\end{tabular}
\caption{Size of the sub-instances of Seville City network.}\label{tab:proportion_OD_pairs}
\end{table}

\subsection{Preliminary experiments}\label{subsec:com_exp_preliminary}

As preliminary results, we show that, for the data used, neither the direct use of \texttt{CPLEX} nor the existing Benders decomposition in \texttt{CPLEX} are competitive with our implementation ad-hoc for the IND problem. 

\texttt{CPLEX} provides three configurations related to its implementation of Benders decomposition. We have set the one that attempts to decompose the model strictly according to the decomposition provided by the user. Table \ref{tab:preliminaries} compares the performance of these three different approaches within a time limit of 4 hours.

\begin{table}[H]
\centering
\begin{tabular}{c|c|c|c|c}
\hline
Procedure & \texttt{Instance$_{g\geq i}$} & \texttt{gap} & \texttt{obj\_v} & \texttt{n\_cuts} \\ \hline
\multirow{5}{*}{\centering \texttt{CPLEX}} 
        & 150 & 170.09 & 1489  & -    \\ 
        & 130 & -      & -     & -    \\ 
        & 100 & -      & -     & -    \\ 
        & 80  & -      & -     & -    \\
        & 60  & -      & -     & -    \\ \hline
\multirow{5}{*}{\centering \texttt{AUTO\_BD}} 
        & 150 & 35.25  & 2377  & 2087 \\ 
        & 130 & 46.13  & 3063  & 2204 \\ 
        & 100 & 101.33 & 3677  & 3807 \\ 
        & 80  & -      & -     & -    \\
        & 60  & -      & -     & -    \\ \hline
\multirow{5}{*}{\centering \texttt{B\&BC}} 
        & 150 & 8.97   & 2522  & 489  \\ 
        & 130 & 13.99  & 3138  & 811  \\ 
        & 100 & 49.04  & 3871  & 1200 \\ 
        & 80 & 57.64 & 4756 & 1763 \\
        & 60  & 141.67 & 5418  & 2240
        \\ \hline
\end{tabular}
\caption{Comparing the performance of \texttt{CPLEX}, \texttt{AUTO\_BD} and \texttt{B\&BC} within a time limit of 4 hours. wE REPORT \texttt{gap} results. The character \texttt{-} denotes those instances in which no solution was found within this time limit.}\label{tab:preliminaries}
\end{table}

It can be observed that not only does \texttt{B\&BC} achieve smaller gaps after 4 hours, but it also finds solutions for the larger instances, which does not happen for the other two procedures. \texttt{B\&BC} did not find solutions within a time limit of 4 hours for bigger instances than the ones shown in this table.

\subsection{Best configuration of Benders decomposition}\label{subsec:comp_exp_benders}

To improve the convergence of our implementation of the branch-and-Benders-cut procedure, and find the best configuration of our Partial Benders decomposition approach described in Subsection \ref{subsec:improved}, we elaborated the following computational experiments related to that decomposition approach.

We performed these experiments by using the already described sub-instances. For each one, we have tested different values of the percentage of O/D pairs to leave in the master problem. We will refer to this parameter as \texttt{Percentage}. That is, we select a percentage of the O/D pairs to leave in the master problem, and with the remaining pairs, we compose the set of Benders sub-problems. The parameter \texttt{Percentage} takes values in the set $\{1,2,5,10,20,40,60\}$. Fixing $\texttt{Percentage}>1$ can be understood as saving useful information in the master problem. We show these results in Table \ref{tab:best_config}, depicting values of gaps. The second column of the table refers to how the percentage of the O/D pairs kept in the master has been selected. If \texttt{Type=1}, we have considered a random selection of them. For \texttt{Type=2}, those O/D pairs with higher demand were selected, and if \texttt{Type=3}, those with lower demand were selected to be kept in the master. We computed the Partial Benders decomposition for each of the sub-instances considered in Table \ref{tab:proportion_OD_pairs}, 90 in total. For each one, the partial Benders decomposition ran for up to 4 hours, in a total of 15 days of computation. Then, in this table, we report the relative optimality gap (\texttt{gap}) obtained after 4 hours. It is shown that for some configurations of the parameters instances \texttt{Instance$_{g\geq 150}$} and \texttt{Instance$_{g\geq 130}$} were solved to optimality using the Partial Benders decomposition, which were not solved using the alternative methods. For these cases, the computational time in seconds is shown in parentheses. Note that for \texttt{Instance$_{g\geq150}$}, the number of O/D pairs is very little, and \texttt{Percentage}=1 does not make sense since this proportion leads to no O/D pair into the master (shown as ``*''). Nevertheless, to get general conclusions, we will examine the details of the largest instances. 

\begin{table}[H]
\footnotesize
    \centering
   \begin{tabular}{c|c|ccccccc}
\hline
\multirow{2}{*}{\texttt{Instance$_{g\geq i}$}} & \multirow{2}{*}{\texttt{Type}} & \multicolumn{6}{c}{\texttt{Percentage}} \\ \cline{3-9} 
\multicolumn{1}{c|}{}   &    & \multicolumn{1}{c}{1}      & \multicolumn{1}{c}{2}      & \multicolumn{1}{c}{5}       & \multicolumn{1}{c}{10}      & \multicolumn{1}{c}{20}      & \multicolumn{1}{c}{40}      &  60 \\ \hline
\multirow{3}{*}{150}                          &\texttt{1}                    & *  & 21.80  & 0 (3870)   & 0 (7091)  & 2.89 & 42.38  & 65.83 \\ 
  &\texttt{2}                    & *   & 0 (9881)   & 14.89   & 0 (11275)  & 0 (14144) & 51.00  & 35.43 \\ 
  &\texttt{3}                    & *  & 0 (4841)   & 0 (5070)   & 29.97  & 39.16  & 49.10  & 54.35 \\ \hline
\multirow{3}{*}{130}                          &\texttt{1}                    & 13.70  & 0 (13385)  & 16.03  & 32.20  & 26.91  & 69.83  & 64.18  \\  
  &\texttt{2}                    & 14.49  & 13.60  & 0 (12345)  & 14.71  & 37.42  & 114.28  & 115.81 \\ 
  &\texttt{3}                    & 32.04  & 0 (12018)  & 0 (8139)  & 10.76  & 56.40  & 82.02  & 88.52  \\ \hline
\multirow{3}{*}{100}                          &\texttt{1}                    & 54.47  & 29.77  & 38.01  & 86.06  & 130.55  & -   & -   \\ 
  &\texttt{2}                    & 52.50  & 29.54  & 32.24  & 60.49  & 145.84  & -   & -   \\ 
  &\texttt{3}                    & 51.80  & 43.53  & 58.83  & 59.05  & 240.90  & -   & -   \\ \hline
\multirow{3}{*}{80}                          &\texttt{1}                    & 83.63   & 66.37  & 90.16  & -  & -   & -   & -   \\ 
  &\texttt{2}                    & 54.42  & 56.44   & 112.21  & 118.11  & -  & -  & -   \\ 
  &\texttt{3}                    & 63.84  & 59.35   & 133.43  & 230.26  & -  & -  & -   \\ \hline
\multirow{3}{*}{60}                          &\texttt{1}                    & 107.18  & 90.88  & -  & -  & -  & -  & -  \\ 
  &\texttt{2}                    & 119.50  & 96.58  & -  & -  & -  & -  & -  \\  
  &\texttt{3}                    & 118.01  & 191.41  & -  & -  & -  & -  & -  \\ \hline
\end{tabular}
    \caption{Comparing the performance of the Partial Benders decomposition for several sub-instances of the Seville City network, with the time limit of 4 hours. We report \texttt{gap} results and also the computational time for those instances that were solved to optimality. With character \texttt{-}, we refer to those instances in which no solution was found keeping this time limit. With character \texttt{*}, we refer to those instances in which it makes no sense this value of \texttt{Percentage} because it corresponds with no selection of O/D pairs.}
    \label{tab:best_config}
\end{table}

By observing the largest instances, it can be seen that fixing 2$\leq\texttt{Percentage}\leq$5 seems to be beneficial, being the configuration \texttt{Percentage}=2 with \texttt{Type}=1 or 2 the best option. Moreover, if \texttt{Percentage}$\geq$10, the problem becomes difficult to solve due to the large size of the master problem. With respect to the parameter \texttt{Type}, it could be that setting it to values 1 or 2 would be more beneficial than in the remaining case due to the following reflection. In these cases, we are considering a certain percentage of those pairs that have a higher demand (\texttt{Type}=2), or part of them that have it (\texttt{Type}=1), to leave them in the master. Since the objective function consists of maximizing the covered demand, and due to the nature of the problem, in the first iteration of \texttt{B\&BC\_P}, for some of these selected O/D pairs, variable $f^w$ is set to 1, which could still be true, although in part, in the optimal solution. This translates into an aid for the convergence of this resolution method.

Finally, you can see that as the number of pairs considered to form the instance increases, the problem becomes more difficult to solve. We have verified that if we consider \texttt{Instance$_{g\geq 40}$}, we will not obtain any solution within 4 hours.

\subsection{Case for the Seville City network}\label{subsec:Comp_exp_case}

\begin{figure}[H]
    \centering
    \begin{subfigure}{0.9\textwidth}
        \centering
        \includegraphics[width=\linewidth]{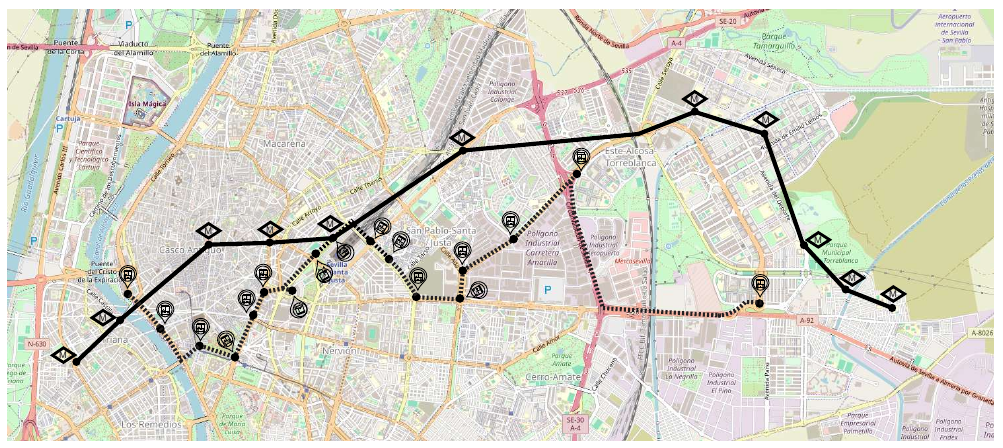}
        \caption{\texttt{Instance}$_{g\geq 80}$}
        \label{fig:subfigura1_optimal}
    \end{subfigure}
    \begin{subfigure}{0.9\textwidth}
        \centering
        \includegraphics[width=\linewidth]{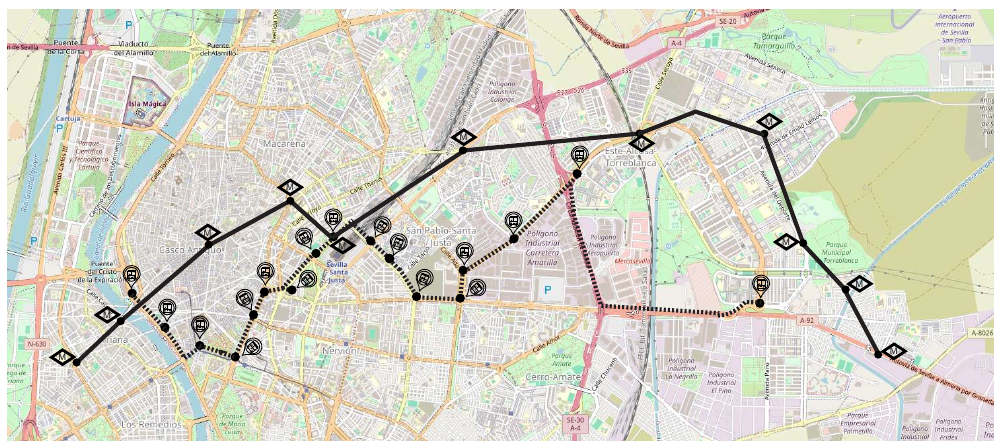}
        \caption{\texttt{Instance}$_{g\geq 60}$}
        \label{fig:subfigura2_optimal}
    \end{subfigure}
    \begin{subfigure}{0.9\textwidth}
        \centering
        \includegraphics[width=\linewidth]{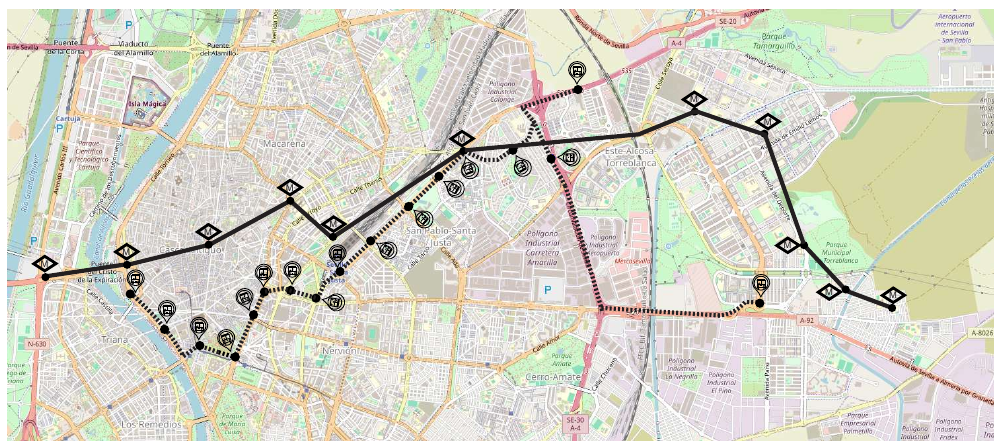}
        \caption{\texttt{Instance}$_{g\geq 40}$}
        \label{fig:subfigura3_optimal}
    \end{subfigure}
    \caption{Optimal network design for some of the sub-instances of the Seville City network.}
    \label{fig:optimal_solutions}
\end{figure}

This sub-section contains the main and final computational experiment. The optimal network designs for some of the sub-instances are shown in Figure \ref{fig:optimal_solutions}. The Figure shows that they are not very different from each other. In fact, the design of mode $\mathcal{R}$ in the largest one, \texttt{Instance}$_{g\geq 40}$, differs only in a few edges from the design in \texttt{Instance}$_{g\geq 80}$ and \texttt{Instance}$_{g\geq 60}$. Concerning the design of mode $\mathcal{S}$, for the two biggest instances considered, it is the same, which mostly coincides with the third one.

Table \ref{tab:optimal_solutions_information} shows computational results and also some insights related to these designs. In terms of computational times, larger instances take a long time to solve to optimality. Parameters \texttt{demand\_R} and \texttt{demand\_S} denote demand covered only by using mode $\mathcal{R}$ and only by using mode $\mathcal{S}$, respectively. The parameter \texttt{demand\_RS} denotes that demand covered that had to use both modes to reach from its origin to its destination. Similarly, parameters \texttt{pairs\_R}, \texttt{pairs\_S} and \texttt{pairs\_RS} denote the number of O/D pairs that only used mode $\mathcal{R}$, those that only used mode $\mathcal{S}$ and those that had to use both modes to be covered, respectively. It can be seen that the majority of the demand covered only needs to use mode $\mathcal{R}$, and that there is more demand that needs to use both modes to be covered than only by using mode $\mathcal{S}$.

\begin{table}[H]
\footnotesize
\centering
\begin{tabular}{c|c|c|c|c|c|c|c|c|c}
\hline
\texttt{Instance$_{g\geq i}$} & \texttt{t} & \texttt{n\_cuts} & \texttt{obj\_v} & \texttt{demand\_R} & \texttt{demand\_S} & \texttt{demand\_RS} & \texttt{pairs\_R} & \texttt{pairs\_S} & \texttt{pairs\_RS} \\ \hline
80              & 45904         & 2553             & 2293                    & 1731               & 0                  & 475                & 14               & 0                & 4                 \\ 
60              & 167885        & 3032             & 3624                    & 3085               & 154                & 385                & 33               & 2                & 4                 \\ 
40              & 443871        & 4446             & 5078                    & 3484               & 204                & 922                & 48               & 3                & 12                \\ \hline
\end{tabular}
\caption{Computational results and some insights related to the optimal network designs for the sub-instances considered in Figure \ref{fig:optimal_solutions}.}
\label{tab:optimal_solutions_information}
\end{table}

%% file: appendix.tex
\begin{appendix}

\section{Linearization of terms \texorpdfstring{$f^{w\mathcal{R}}_a z^{\mathcal{R}}_k$}{} and \texorpdfstring{$f^{w\mathcal{S}}_a z^{\mathcal{S}}_k$}{}}\label{append:linearization}

We define new variables $h^{w\mathcal{R}}_a = f^{w\mathcal{R}}_a z^{\mathcal{R}}_k$ and $h^{w\mathcal{S}}_a = f^{w\mathcal{S}}_a z^{\mathcal{S}}_k$.
    \begin{flalign}
    & h^{w\mathcal{R}}_a \leq f^{w\mathcal{R}}_a, \quad w\in W, a\in A_{\mathcal{R}},&\label{eq_apendixA:linearizationR_1}\\
    & h^{w\mathcal{R}}_a \leq z^{\mathcal{R}}_{a^s}, \quad w\in W, a\in A_{\mathcal{R}},&\label{eq_apendixA:linearizationR_2}\\
    & f^{w\mathcal{R}}_a - \left(1 - z^{\mathcal{R}}_{a^s}\right) \leq h^{w\mathcal{R}}_a, \quad w\in W, a\in A_{\mathcal{R}},&\label{eq_apendixA:linearizationR_3}\\
    & h^{w\mathcal{R}}_a \in \{0,1\}, w\in W, a\in A_{\mathcal{R}}.&
    \end{flalign}
    \begin{flalign}
    & h^{w\mathcal{S}}_a \leq f^{w\mathcal{S}}_a, \quad w\in W, a\in A_{\mathcal{S}},&\label{eq_apendixA:linearizationS_1}\\
    & h^{w\mathcal{S}}_a \leq z^{\mathcal{S}}_{a^s}, \quad w\in W, a\in A_{\mathcal{S}},&\label{eq_apendixA:linearizationS_2}\\
    & f^{w\mathcal{S}}_a - \left(1 - z^{\mathcal{S}}_{a^s}\right) \leq h^{w\mathcal{S}}_a, \quad w\in W, a\in A_{\mathcal{S}},&\label{eq_apendixA:linearizationS_3}\\
    & h^{w\mathcal{S}}_a \in \{0,1\}, w\in W, a\in A_{\mathcal{S}}.&
    \end{flalign}

    


\section{Benders formulation}\label{append:benders}

In the following, we describe a Benders implementation obtained by projecting out variables $f^{w\mathcal{R}}_a$, $f^{w\mathcal{S}}_a$, $f^{w\mathcal{SR}}_k$ and $f^{w\mathcal{RS}}_k$.

The master problem $(MP)$ that we solve is:
\begin{align}
(MP) \quad \max\limits_{\boldsymbol{x},\boldsymbol{y},\boldsymbol{z},\boldsymbol{f},\boldsymbol{v},\boldsymbol{h}} &\sum_{w\in W} g^w f^w& \\
\mbox{s.t. }& \mbox{(\ref{eq:budget})-(\ref{eq:relation_8})},(\ref{eq_apendixA:linearizationR_2}),(\ref{eq_apendixA:linearizationS_2}) & \nonumber\\
&+\{\mbox{Benders Cuts }(\boldsymbol{x},\boldsymbol{y},\boldsymbol{z}, \boldsymbol{f})\}& \nonumber\\
&x^{\mathcal{R}}_e,\, x_e^{\mathcal{S}},\, y^{\mathcal{R}}_k,\, z^{\mathcal{R}}_k,\, z_k^{\mathcal{S}},\, f^w, h^{w\mathcal{R}}_a, h^{w\mathcal{S}}_a,\, v_{{w^s}k}^{w\mathcal{R}}, \,v_{{w^s}k}^{w\mathcal{S}}, \,v_{k{w^t}}^{w\mathcal{R}}, \,v_{k{w^t}}^{w\mathcal{S}} \in \{0,1\}. & \nonumber
\end{align}

Since the structure of the model allows it, we consider a feasibility subproblem made of constraints \eqref{eq:relation_3}-\eqref{eq:utility},\eqref{eq_apendixA:linearizationR_1},\eqref{eq_apendixA:linearizationR_3},\eqref{eq_apendixA:linearizationS_1},\eqref{eq_apendixA:linearizationS_3}, for each commodity $w\in W$, and denoted by $(SP)^w$. Then, we must consider the following subproblem structure.

\begin{itemize}
    \item Flow conservation constraints
    \begin{flalign}
&\sum_{a \in \delta_w^+(k)} f^{\mathcal{R}}_a - \sum_{a \in \delta_w^-(k)} f^{\mathcal{R}}_a =
v_{{w^s}k}^{\mathcal{R},out} - v_{k{w^t}}^{\mathcal{R},out} - \lambda \left(\Delta\, v_{{w^s}k}^{\mathcal{R}} - \Delta\, v_{k{w^t}}^{\mathcal{R}}\right), \quad k\in N_{\mathcal{R}} \setminus N_{trans}, &
    \label{eq:flow1}\\
&\sum_{a \in \vartheta_w^+(k)} f^{\mathcal{S}}_a - \sum_{a \in \vartheta_w^-(k)} f^{\mathcal{S}}_a =
v_{{w^s}k}^{\mathcal{S},out} - v_{k{w^t}}^{\mathcal{S},out} - \lambda \left(\Delta\, v_{{w^s}k}^{\mathcal{S}} - \Delta\, v_{k{w^t}}^{\mathcal{S}}\right), \quad k\in N_{\mathcal{S}} \setminus N_{trans}, &
    \label{eq:flow11}\\
					& v_{{w^s}k}^{\mathcal{R},out}-\lambda\,\Delta\, v_{{w^s}k}^{\mathcal{R}} \leq \sum_{a \in \delta_w^+(k)} f^{w\mathcal{R}}_a,\quad w\in W, \, k\in N_{\mathcal{R}},& \label{eq:flow2}\\
					& v_{{w^s}k}^{\mathcal{S},out} -\lambda\,\Delta\,v_{{w^s}k}^{\mathcal{S}}\leq \sum_{a \in \vartheta_w^+(k)} f^{w\mathcal{S}}_a,\quad w\in W, \, k\in N_{\mathcal{S}}, &\label{eq:flow3}\\
					& v_{k{w^t}}^{\mathcal{R},out} -\lambda\,\Delta\,v_{k{w^t}}^{\mathcal{R}}\leq \sum_{a \in \delta_w^-(k)} f^{w\mathcal{R}}_a,\quad w\in W, \, k\in N_{\mathcal{R}}, &\label{eq:flow4}\\
					& v_{k{w^t}}^{\mathcal{S},out} -\lambda\,\Delta\,v_{k{w^t}}^{\mathcal{S}}\leq \sum_{a \in \vartheta_w^-(k)} f^{w\mathcal{S}}_a,\quad w\in W, \, k\in N_{\mathcal{S}}, &\label{eq:flow5}\\
         &f^{\mathcal{R}}_a\leq f^{out}-\lambda\,\Delta\,f, \quad a\in A_{\mathcal{R}}, &\label{eq:flow6}\\
    &f^{\mathcal{S}}_a\leq f^{out}-\lambda\,\Delta\,f, \quad a\in A_{\mathcal{S}}.& \label{eq:flow7}
    \end{flalign}
    \item Transfer constraints
    \begin{flalign}
    &\sum_{k \in N_{trans}} f^{\mathcal{SR}}_k \leq 1,  &\label{eq:trans}\\
    &\sum_{k \in N_{trans}} f^{\mathcal{RS}}_k \leq 1,  &\label{eq:trans2}\\
    & \sum_{a \in \delta_w^-(k)} f^{\mathcal{R}}_a + f^{\mathcal{SR}}_k - \sum_{a \in \delta_w^+(k)} f^{\mathcal{R}}_a = -\left(v_{{w^s}k}^{\mathcal{R},out}-\lambda\,\Delta\,v_{{w^s}k}^{\mathcal{R}}\right)+\left(v_{k{w^t}}^{\mathcal{R},out}-\lambda\,\Delta\,v_{k{w^t}}^{\mathcal{R}}\right), \quad k\in N_{trans},& \label{eq:trans3}\\
    &\sum_{a \in \vartheta_w^-(k)} f^{\mathcal{S}}_a + f^{\mathcal{RS}}_k - \sum_{a \in \vartheta_w^+(k)} f^{\mathcal{S}}_a  = -\left(v_{{w^s}k}^{\mathcal{S},out}-\lambda\,\Delta\,v_{{w^s}k}^{\mathcal{S}}\right)+\left(v_{k{w^t}}^{\mathcal{S},out}-\lambda\,\Delta\,v_{k{w^t}}^{\mathcal{S}}\right), \quad k\in N_{trans}.&\label{eq:trans4}
    \end{flalign}
    \begin{flalign}
& f^{w\mathcal{\mathcal{RS}}}_k \leq \sum_{a \in \vartheta_w^+(k)} f^{w\mathcal{S}}_a, \qquad f^{w\mathcal{\mathcal{RS}}}_k \leq \sum_{a \in \delta_w^-(k)} f^{w\mathcal{R}}_a, \quad w\in W, \, k\in N_{trans}, \label{eq:trans5} &\\
& f^{w\mathcal{\mathcal{SR}}}_k \leq \sum_{a \in \vartheta_w^-(k)} f^{w\mathcal{S}}_a, \qquad f^{w\mathcal{\mathcal{SR}}}_k \leq \sum_{a \in \delta_w^+(k)} f^{w\mathcal{R}}_a, \quad w\in W, \, k\in N_{trans}. \label{eq:trans6} &
 \end{flalign}
    \item Location-allocation constraints
    \begin{flalign}
    &f^{\mathcal{R}}_a + f^{\mathcal{R}}_{\hat{a}} \leq x^{\mathcal{R},out}_e - \lambda\,\Delta\,x^{\mathcal{R}}_e, \quad e= \{i,j\}\in E_{\mathcal{R}}: a=(i,j), \hat{a}=(j,i),& \label{eq:loc}\\
    &f^{\mathcal{S}}_a + f^{\mathcal{S}}_{\hat{a}} \leq x^{\mathcal{S},out}_e-\lambda\,\Delta\,x^{\mathcal{S}}_e, \quad e= \{i,j\}\in E_{\mathcal{S}}: a=(i,j), \hat{a}=(j,i),& \label{eq:loc2}\\
     &f^{w\mathcal{R}}_a + f^{w\mathcal{R}}_{\hat{a}}+f^{w\mathcal{S}}_a + f^{w\mathcal{S}}_{\hat{a}} \leq 1, \quad w\in W, e= \{i,j\}\in E_{\mathcal{R}}\cap E_{\mathcal{S}}: a=(i,j), \hat{a}=(j,i).&\label{eq:loc3}
    \end{flalign}
    \item Alignment stop constraints
    \begin{flalign}
    &f^{\mathcal{SR}}_k + f^{\mathcal{RS}}_k\leq z^{\mathcal{R},out}_k-\lambda\Delta\,z^{\mathcal{R}}_k, \quad k\in N_{trans},& \label{eq:align}\\
    &f^{\mathcal{SR}}_k + f^{\mathcal{RS}}_k\leq z^{\mathcal{S},out}_k -\lambda\Delta\,z^{\mathcal{S}}_k, \quad k\in N_{trans}.& \label{eq:align2}
    \end{flalign}
    \item Mode choice constraint
    \begin{flalign}
    \begin{split}
    &\sum_{k \in N_{\mathcal{R}}} \left( t_{{w^s}k}\left( v_{{w^s}k}^{\mathcal{R},out}-\lambda\,\Delta\, v_{{w^s}k}^{\mathcal{R}}\right) + t_{k{w^t}}\left(v_{k{w^t}}^{\mathcal{R},out} -\lambda\,\Delta\,v_{k{w^t}}^{\mathcal{R}}\right) \right)+ \\
    &+\sum_{k \in N_{\mathcal{S}}} \left( t_{{w^s}k}\left( v_{{w^s}k}^{\mathcal{S},out}-\lambda\,\Delta\, v_{{w^s}k}^{\mathcal{S}}\right) + t_{k{w^t}}\left(v_{k{w^t}}^{\mathcal{S},out} -\lambda\,\Delta\,v_{k{w^t}}^{\mathcal{S}}\right) \right)+ \\
    &+\sum_{a\in A_{\mathcal{R}}}t^{\mathcal{R}}_af^{\mathcal{R}}_a + \sum_{a\in A_{\mathcal{S}}}t^{\mathcal{S}}_af^{\mathcal{S}}_a + \sum_{k\in N_{trans}}t^{\mathcal{RS}}_kf^{\mathcal{RS}}_k + \sum_{k\in N_{trans}}t^{\mathcal{SR}}_kf^{\mathcal{SR}}_k + \\
    & + t^{\mathcal{R}}_{stop}\sum_{k\in N_{\mathcal{R}}}\sum_{a\in\delta^+(k)}\left(h_a^{\mathcal{R},out}-\lambda\,\Delta\,h_a^{\mathcal{R}}\right) + t^{\mathcal{S}}_{stop}\sum_{k\in N_{\mathcal{S}}}\sum_{a\in\gamma^+(k)}\left(h_a^{\mathcal{S},out}-\lambda\,\Delta\,h_a^{\mathcal{S}}\right) +\\ &+\left(t^{\mathcal{R}}_{wait}-\frac{1}{2}t^{\mathcal{R}}_{stop}\right)\left(f^{out}-\lambda\,\Delta\,f\right)\leq u_{priv}\left(f^{out}-\lambda\,\Delta\,f\right),\label{eq:uti}
    \end{split}\end{flalign}
    
    \item Linearization of terms $f^{\mathcal{R}}_a z^{\mathcal{R}}_k$ and $f^{\mathcal{S}}_a z^{\mathcal{S}}_k$.
    
    We define new variables $h^{\mathcal{R}}_a = f^{\mathcal{R}}_a z^{\mathcal{R}}_k$ and $h^{\mathcal{S}}_a = f^{\mathcal{S}}_a z^{\mathcal{S}}_k$.
    \begin{flalign}
    & h^{\mathcal{R},out}_a - \lambda\,\Delta h^{\mathcal{R}}_a \leq f^{\mathcal{R}}_a, \quad a\in A_{\mathcal{R}},&\\
    & h^{\mathcal{R},out}_a - \lambda\,\Delta h^{\mathcal{R}}_a\leq z^{\mathcal{R},out}_{a^s} - \lambda\,\Delta z^{\mathcal{R}}_{a^s}, \quad a\in A_{\mathcal{R}},&\\
    & f^{\mathcal{R}}_a - \left(1 - \left(z^{\mathcal{R},out}_{a^s} - \lambda\,\Delta z^{\mathcal{R}}_{a^s}\right)\right) \leq h^{\mathcal{R},out}_a - \lambda\,\Delta h^{\mathcal{R}}_{a}, \quad a\in A_{\mathcal{R}},&\\
    & h^{\mathcal{R}}_a \in \{0,1\}, a\in A_{\mathcal{R}}.
    \end{flalign}
    \begin{flalign}
    & h^{\mathcal{S},out}_a - \lambda\,\Delta h^{\mathcal{S}}_a \leq f^{\mathcal{S}}_a, \quad a\in A_{\mathcal{S}},&\\
    & h^{\mathcal{S},out}_a - \lambda\,\Delta h^{\mathcal{S}}_a \leq z^{\mathcal{S},out}_{a^s} - \lambda\,\Delta z^{\mathcal{S}}_{a^s}, \quad a\in A_{\mathcal{S}},&\\
    & f^{\mathcal{S}}_a - \left(1 - \left(z^{\mathcal{S},out}_{a^s} - \lambda\,\Delta z^{\mathcal{S}}_{a^s}\right)\right) \leq h^{\mathcal{S},out}_a - \lambda\,\Delta h^{\mathcal{S}}_{a}, \quad a\in A_{\mathcal{S}},&\\
    & h^{\mathcal{S}}_a \in \{0,1\}, a\in A_{\mathcal{S}}.
    \label{eq:linearization_S}
    \end{flalign}

    \item Continuous variables
    \begin{flalign}
    & f_a^{\mathcal{R}},\, f_a^{\mathcal{S}},\, f_k^{\mathcal{RS}},\, f_k^{\mathcal{SR}} \in [0,1]&
    \end{flalign}
    
\end{itemize}


The dual of each feasibility subproblem can be expressed as follows, where $\boldsymbol{\alpha}^{\mathcal{R}}$ and $\boldsymbol{\alpha}^{\mathcal{S}}$ are the vectors of dual variables of constraints \eqref{eq:flow1} and \eqref{eq:flow11}, $\boldsymbol{\epsilon}^{\mathcal{R}}$ and $\boldsymbol{\epsilon}^{\mathcal{S}}$ are vectors of dual variables corresponding to the set of constraints \eqref{eq:flow2}-\eqref{eq:flow5}, $\boldsymbol{\delta}^{\mathcal{R}}$ and $\boldsymbol{\delta}^{\mathcal{S}}$ are the vectors of dual variables of constraints \eqref{eq:flow6} and \eqref{eq:flow7}. The dual variables vectors $\boldsymbol{\nu}^{\mathcal{RS}}$ and $\boldsymbol{\nu}^{\mathcal{SR}}$ correspond to constraints \eqref{eq:trans} and \eqref{eq:trans2}. The sets of constraints \eqref{eq:trans3} and \eqref{eq:trans4} are identified with the dual variables vectors $\boldsymbol{\beta}^{\mathcal{R}}$ and $\boldsymbol{\beta}^{\mathcal{S}}$. The dual variables vectors $\boldsymbol{\phi}^{\mathcal{RS}}$ and $\boldsymbol{\phi}^{\mathcal{SR}}$ correspond to constraints \eqref{eq:trans5} and \eqref{eq:trans6}. The dual vectors $\boldsymbol{\sigma}^{\mathcal{R}}$, $\boldsymbol{\sigma}^{\mathcal{S}}$ and $\boldsymbol{\sigma}^{\mathcal{RS}}$ are related with constraints \eqref{eq:loc}, \eqref{eq:loc2} and \eqref{eq:loc3}. The set of constraints \eqref{eq:align} and \eqref{eq:align2} are identified with the dual variables vectors $\boldsymbol{\theta}^{\mathcal{R}}$ and $\boldsymbol{\theta}^{\mathcal{S}}$, respectively. Finally, vector $\boldsymbol{\eta}$ represents the set of constraints \eqref{eq:uti}. Furthermore, the dual variables vectors $\boldsymbol{\lambda}$ and $\boldsymbol{\rho}$ correspond with constraints related to the linearization part of the mode choice constraint.

\begin{itemize}
    \item Objective function:
    \begin{flalign}
    \begin{split}
&\max_{\boldsymbol{\alpha},\boldsymbol{\beta},\boldsymbol{\gamma},\boldsymbol{\delta},\boldsymbol{\sigma},\boldsymbol{\eta},\boldsymbol{\nu},\boldsymbol{\theta}}\sum_{k\in N_\mathcal{R}\setminus N_{trans}}\left(v_{{w^s}k}^{\mathcal{R},out} - v_{k{w^t}}^{\mathcal{R},out} \right)\,\alpha_{k}^{\mathcal{R}} + \sum_{k\in N_\mathcal{S}\setminus N_{trans}}\left(v_{{w^s}k}^{\mathcal{S},out} - v_{k{w^t}}^{\mathcal{S},out}\right)\,\alpha_{k}^{\mathcal{S}} +\\
    & + \sum_{k\in N_\mathcal{R}} v_{{w^s}k}^{\mathcal{R},out}\,\epsilon_{k}^{\mathcal{R}} + \sum_{k\in N_\mathcal{S}} v_{{w^s}k}^{\mathcal{S},out}\,\epsilon_{k}^{\mathcal{S}} + \sum_{k\in N_\mathcal{R}} v_{k{w^t}}^{\mathcal{R},out}\,\epsilon_{k}^{\mathcal{R}} + \sum_{k\in N_\mathcal{S}} v_{k{w^t}}^{\mathcal{S},out}\,\epsilon_{k}^{\mathcal{S}}+ \\
    & - \sum_{e\in E_{\mathcal{R}}\cap E_{\mathcal{S}}}\sigma_e^{\mathcal{RS}} - \sum_{e\in E_{\mathcal{R}}}x_e^{\mathcal{R},out}\,\sigma_e^{\mathcal{R}} - \sum_{e\in E_{\mathcal{S}}}x_e^{\mathcal{S},out}\,\sigma_e^{\mathcal{S}} -\nu^{\mathcal{SR}} - \nu^{\mathcal{RS}} - \\ 
    & - \sum_{a\in A_{\mathcal{R}}}f^{out}\,\delta_a^{\mathcal{R}} - \sum_{a\in A_{\mathcal{S}}}f^{out}\,\delta_a^{\mathcal{S}} - \sum_{k\in N_{trans}}z_k^{\mathcal{R},out}\,\theta_k^{\mathcal{R}} - \sum_{k\in N_{trans}}z_k^{\mathcal{S},out}\,\theta_k^{\mathcal{S}} - \\
    & - \left(f^{out}\,\left(u -t^{\mathcal{R}}_{wait} + \frac{1}{2}\,t^{\mathcal{R}}_{stop}\right) - t^{\mathcal{R}}_{stop}\sum_{a\in A_{\mathcal{R}}}h^{\mathcal{R},out}_a - t^{\mathcal{S}}_{stop}\sum_{a\in A_{\mathcal{S}}}h^{\mathcal{S},out}_a+ \right.\\
    &\left. -\sum_{k \in N_{\mathcal{R}}} \left( t_{{w^s}k}v_{{w^s}k}^{\mathcal{R},out} + t_{k{w^t}}v_{k{w^t}}^{\mathcal{R},out} \right) - \sum_{k \in N_{\mathcal{S}}} \left( t_{{w^s}k}v_{{w^s}k}^{\mathcal{S},out} + t_{k{w^t}}v_{k{w^t}}^{\mathcal{S},out} \right)\right)\eta \, +\\
    & + \sum_{k\in N_{trans}}\left(-\left(v_{{w^s}k}^{\mathcal{R},out} - v_{k{w^t}}^{\mathcal{R},out} \right)\beta_k^{\mathcal{R}} + \left(v_{{w^s}k}^{\mathcal{S},out}- v_{k{w^t}}^{\mathcal{S},out}\right)\beta_k^{\mathcal{S}}\right) + \\
    & + \sum_{a \in A_{\mathcal{R}}} h^{\mathcal{R},out}_a\,\lambda^{\mathcal{R}}_a + \sum_{a \in A_{\mathcal{S}}} h^{\mathcal{S},out}_a\,\lambda^{\mathcal{S}}_a-  \\
    &  - \sum_{a\in A_{\mathcal{R}}}\left(1 + h^{\mathcal{R},out}_a - z^{\mathcal{R},out}_{a^s}\right)\,\rho^{\mathcal{R}}_a - \sum_{a\in A_{\mathcal{S}}}\left(1 + h^{\mathcal{S},out}_a - z^{\mathcal{S},out}_{a^s}\right)\,\rho^{\mathcal{S}}_a
    \end{split}
    \end{flalign}
    
    \item Normalization constraint
    \begin{flalign}
    \begin{split}
    & \sum\limits_{j\in N_{\mathcal{R}} \setminus N_{trans}}\left(\Delta\, v_{{w^s}j}^{\mathcal{R}}- \Delta\, v_{j{w^t}}^{\mathcal{R}} -\right)\alpha_{j}^{\mathcal{R}} +\sum\limits_{j\in N_{\mathcal{S}}\setminus N_{trans}}\left(\Delta\, v_{{w^s}j}^{\mathcal{S}} - \Delta\, v_{j{w^t}}^{\mathcal{S}}\right)\alpha_{j}^{\mathcal{S}} +\\
    & +\sum\limits_{j\in N_{\mathcal{R}}}\left(\Delta\, v_{{w^s}j}^{\mathcal{R}}+\Delta\, v_{j{w^t}}^{\mathcal{R}}\right)\epsilon_j^{\mathcal{R}}+\sum\limits_{j\in N_{\mathcal{S}}}\left(\Delta\, v_{{w^s}j}^{\mathcal{S}}+\Delta\, v_{j{w^t}}^{\mathcal{S}}\right)\epsilon_j^{\mathcal{S}}-\\
    &- \sum_{e\in E_{\mathcal{R}}}\Delta\,x^{\mathcal{R}}_e\sigma^{\mathcal{R}}_e - \sum_{e\in E_{\mathcal{S}}}\Delta\,x^{\mathcal{S}}_e\sigma^{\mathcal{S}}_e - \sum_{a\in A_{\mathcal{R}}}\Delta\,f\,\delta^{\mathcal{R}}_a - \sum_{a\in A_{\mathcal{S}}}\Delta\,f\delta^{\mathcal{R}}_a - \\
    & - \sum_{k\in N_{trans}}\Delta\,z^{\mathcal{R}}_k\theta^{\mathcal{R}}_k - \sum_{k\in N_{trans}}\Delta\,z^{\mathcal{S}}_k\theta^{\mathcal{S}}_k + \\
    & + \sum_{k\in N_{trans}}\left(\left(\Delta v_{{w^s}k}^{\mathcal{R}} - \Delta v_{k{w^t}}^{\mathcal{R}} \right)\beta_k^{\mathcal{R}} + \left(\Delta v_{{w^s}k}^{\mathcal{S}}- \Delta v_{k{w^t}}^{\mathcal{S}}\right)\beta_k^{\mathcal{S}}\right)+ \\
    & + \left(\left(-u+t^{\mathcal{R}}_{wait}-\frac{1}{2}t^{\mathcal{R}}_{stop}\right)\Delta\,f + t^{\mathcal{R}}_{stop}\sum_{k\in N_{\mathcal{R}}}\sum_{a\in \delta^+(k)}\Delta\,h^{\mathcal{R}}_{a} + t^{\mathcal{S}}_{stop}\sum_{k\in N_{\mathcal{S}}}\sum_{a\in \vartheta^+(k)}\Delta\,h^{\mathcal{S}}_{a}\right.+\\
    & +\left.\sum_{k \in N_{\mathcal{R}}} \left( t_{{w^s}k}\Delta\, v_{{w^s}k}^{\mathcal{R}} + t_{k{w^t}}\Delta\,v_{k{w^t}}^{\mathcal{R}} \right)+
    \sum_{k \in N_{\mathcal{S}}} \left( t_{{w^s}k}\Delta\, v_{{w^s}k}^{\mathcal{S}} + t_{k{w^t}}\Delta\,v_{k{w^t}}^{\mathcal{S}} \right) \right)\eta\, +\\
    & + \sum_{a\in A_{\mathcal{R}}}\Delta h^{\mathcal{R}}_a \lambda^{\mathcal{R}}_a + \sum_{a\in A_{\mathcal{S}}}\Delta h^{\mathcal{S}}_a \lambda^{\mathcal{S}}_a +\\ 
    & 
    + \sum_{a\in A_{\mathcal{R}}}\left(\Delta z^{\mathcal{R}}_{a^s} - \Delta h^{\mathcal{R}}_a\right)\rho^{\mathcal{R}}_a + \sum_{a\in A_{\mathcal{S}}}\left(\Delta z^{\mathcal{S}}_{a^s} - \Delta h^{\mathcal{S}}_a\right)\rho^{\mathcal{S}}_a \leq 1
    \end{split}
    \end{flalign}
\end{itemize}

\begin{itemize}
    \item Constraints:
        \begin{equation}
        \left.\begin{aligned}
        &\epsilon_i^{\mathcal{R}} + \epsilon_j^{\mathcal{R}} - \sigma^{\mathcal{R}}_e - \delta^{\mathcal{R}}_a - t^{\mathcal{R}}_a\,\eta + \lambda^{\mathcal{R}}_a - \rho^{\mathcal{R}}_a \\
        &\circ \text{ if } i\notin N_{trans}: \quad +\,\alpha_i^{\mathcal{R}}\\
        &\circ \text{ if } j\notin N_{trans}: \quad -\,\alpha_j^{\mathcal{R}}\\
        &\circ \text{ if } j\in N_{trans}: \quad + \,\beta^{\mathcal{R}}_j\\
        &\circ \text{ if } i\in N_{trans}: \quad - \,\beta^{\mathcal{R}}_i\\
        &\circ \text{ if } e\in E_{\mathcal{R}}\cap E_{\mathcal{S}}: \quad -\sigma_e^{\mathcal{SR}}\\
        &\circ \text{ if } i\in N_{trans}: \quad + \, \phi_i^{\mathcal{RS}}\\
        &\circ \text{ if } j\in N_{trans}: \quad + \, \phi_j^{\mathcal{SR}}
        \end{aligned}\right\}\leq 0, \quad a=(i,j)\in A_{\mathcal{R}}:\, e=\{i,j\},
        \end{equation}
        
        \begin{equation}
        \left.\begin{aligned}
        & \epsilon_i^{\mathcal{S}} + \epsilon_j^{\mathcal{S}} - \sigma^{\mathcal{S}}_e - \delta^{\mathcal{S}}_a - t^{\mathcal{S}}_a\,\eta + \lambda^{\mathcal{S}}_a - \rho^{\mathcal{S}}_a\\
        &\circ \text{ if } i\notin N_{trans}: \quad +\,\alpha_i^{\mathcal{S}}\\
        &\circ \text{ if } j\notin N_{trans}: \quad -\,\alpha_j^{\mathcal{S}}\\
        &\circ \text{ if } j\in N_{trans}: \quad + \,\beta^{\mathcal{S}}_j\\
        &\circ \text{ if } i\in N_{trans}: \quad - \,\beta^{\mathcal{S}}_i\\
        &\circ \text{ if } e\in E_{\mathcal{R}}\cap E_{\mathcal{S}}: \quad -\sigma_e^{\mathcal{SR}}\\
        &\circ \text{ if } i\in N_{trans}: \quad + \, \phi_i^{\mathcal{RS}}\\
        &\circ \text{ if } j\in N_{trans}: \quad + \, \phi_j^{\mathcal{SR}}
        \end{aligned}\right\}\leq 0, \quad a=(i,j)\in A_{\mathcal{S}}:\, e=\{i,j\},
        \end{equation}

        \begin{equation}
        -\nu^{\mathcal{SR}} + \beta^{\mathcal{R}}_k - \beta^{\mathcal{S}}_k - \phi_k^{\mathcal{RS}} - \phi_k^{\mathcal{RS}} - \theta^{\mathcal{R}}_k - \theta^{\mathcal{S}}_k - t^{\mathcal{SR}}_k\,\eta \leq 0, \quad k\in N_{trans},
        \end{equation}
        \begin{equation}
        -\nu^{\mathcal{RS}} + \beta^{\mathcal{S}}_k - \beta^{\mathcal{R}}_k - \phi_k^{\mathcal{SR}} - \phi_k^{\mathcal{SR}} -\theta^{\mathcal{R}}_k - \theta^{\mathcal{S}}_k - t^{\mathcal{RS}}_k\,\eta \leq 0, \quad k\in N_{trans},
        \end{equation}
        
        \begin{equation}
        \sigma^{\mathcal{R}}_e,\, \sigma^{\mathcal{S}}_e,\,
            \sigma^{\mathcal{RS}}_e,\,\phi_i^{\mathcal{RS}},\, \phi_i^{\mathcal{SR}},\,
            \epsilon^{\mathcal{R}}_i,\, 
            \epsilon^{\mathcal{S}}_i,\, \delta^{\mathcal{R}}_a,\, \delta^{\mathcal{S}}_a,\, \theta^{\mathcal{R}}_i,\, \theta^{\mathcal{S}}_i,\, \nu^{\mathcal{RS}},\, \nu^{\mathcal{SR}},\, \eta \geq 0.
        \end{equation}
\end{itemize}


Benders cuts:

\begin{equation}
\left.\begin{aligned}
& \left( - \sum_{a\in A_{\mathcal{R}}}\delta_a^{\mathcal{R}} - \sum_{a\in A_{\mathcal{S}}}\delta_a^{\mathcal{S}} + \left(-u+t^{\mathcal{R}}_{wait} - \frac{1}{2}\,t^{\mathcal{R}}_{stop}\right)\eta \right)\,f +\\
& + \sum_{k \in N_{\mathcal{R}}\setminus N_{trans}} \left( \left(\alpha_k^{\mathcal{R}} +\eta\, t_{{w^s}k}+\epsilon_{{w^s}k}^{\mathcal{R}}\right)v_{{w^s}k}^{\mathcal{R}} + \left(- \alpha_k^{\mathcal{R}} + \eta\, t_{k{w^t}}+\epsilon_{k{w^t}}^{\mathcal{R}}\right)v_{k{w^t}}^{\mathcal{R}} \right)+\\
& + \sum_{k \in N_{\mathcal{S}}\setminus N_{trans}} \left( \left( \alpha_k^{\mathcal{S}} + \eta\, t_{{w^s}k}+\epsilon_{{w^s}k}^{\mathcal{S}}\right)v_{{w^s}k}^{\mathcal{S}} + \left(- \alpha_k^{\mathcal{S}} +\eta\, t_{k{w^t}}+\epsilon_{k{w^t}}^{\mathcal{S}}\right)v_{k{w^t}}^{\mathcal{S}} \right) +\\
& + \sum_{k \in N_{trans}}  \left((\eta\, t_{{w^s}k}+\epsilon_{{w^s}k}^{\mathcal{R}}-\beta_k^{\mathcal{R}})\,v_{{w^s}k}^{\mathcal{R}} + (\eta\, t_{k{w^t}}+\epsilon_{k{w^t}}^{\mathcal{R}}+\beta_k^{\mathcal{R}})\,v_{k{w^t}}^{\mathcal{R}} \right) + \\
& + \sum_{k \in N_{trans}}  \left((\eta\, t_{{w^s}k}+\epsilon_{{w^s}k}^{\mathcal{S}}-\beta_k^{\mathcal{S}})\,v_{{w^s}k}^{\mathcal{S}} + (\eta\, t_{k{w^t}}+\epsilon_{k{w^t}}^{\mathcal{S}}+\beta_k^{\mathcal{S}})\,v_{k{w^t}}^{\mathcal{S}} \right) + \\
& - \sum_{e\in E_{\mathcal{R}}}\sigma_e^{\mathcal{R}} \, x_e^{\mathcal{R}}- \sum_{e\in E_{\mathcal{S}}}\sigma_e^{\mathcal{S}} \, x_e^{\mathcal{S}} +\\
& + \sum_{a\in A_{\mathcal{R}}}\left(\lambda^{\mathcal{R}}_a + t^{\mathcal{R}}_{stop}\eta - \rho^{\mathcal{R}}_a\right)\, h^{\mathcal{R}}_a + \sum_{a\in A_{\mathcal{S}}}\left(\lambda^{\mathcal{S}}_a + t^{\mathcal{S}}_{stop}\eta - \rho^{\mathcal{S}}_a\right)\,h^{\mathcal{S}}_a +\\
& + \sum_{k\in N_{trans}} \left(- \theta^{\mathcal{R}}_k + \sum_{a\in \delta^+(k)}\rho^{\mathcal{R}}_a\right)z^{\mathcal{R}}_k + \sum_{k\in N_{\mathcal{R}}\setminus \{N_{trans}\}} \left( \sum_{a\in \delta^+(k)}\rho^{\mathcal{R}}_a \right)z^{\mathcal{R}}_k +\\
& + \sum_{k\in N_{trans}} \left(- \theta^{\mathcal{S}}_k + \sum_{a\in \vartheta^+(k)}\rho^{\mathcal{S}}_a \right)z^{\mathcal{S}}_k+ \sum_{k\in N_{\mathcal{S}}\setminus \{N_{trans}\}} \left( \sum_{a\in \vartheta^+(k)}\rho^{\mathcal{S}}_a \right)z^{\mathcal{S}}_k
\end{aligned}\right\}\leq \ast\ast
\end{equation}

\begin{equation*}
    \ast\ast = \nu^{\mathcal{SR}} + \nu^{\mathcal{RS}} + \sum_{a\in A_{\mathcal{R}}} \rho^{\mathcal{R}}_a + \sum_{a\in A_{\mathcal{S}}} \rho^{\mathcal{S}}_a + \sum_{e\in E_{\mathcal{R}}\cap E_{\mathcal{S}}}\sigma_e^{\mathcal{RS}}
\end{equation*}

\end{appendix}